\documentclass[a4paper,12pt,reqno]{amsart}
\usepackage{amsmath}
\usepackage{amssymb}
\usepackage{amsthm}

\usepackage{amscd}
\usepackage{amsfonts}
\usepackage{setspace}
\usepackage{version}

\title[Minor arcs and mean values for smooth numbers]{Minor arcs, mean values, and restriction theory for exponential sums over smooth numbers}
\author{Adam J Harper}
\address{Jesus College, Cambridge CB5 8BL, England}
\email{A.J.Harper@dpmms.cam.ac.uk}
\date{7th August 2014}
\thanks{The author is supported by a research fellowship at Jesus College, Cambridge.}

\numberwithin{equation}{section}
\theoremstyle{plain}

\onehalfspace
\newcommand{\N}{\mathbb{N}}
\newcommand{\R}{\mathbb{R}}

\newcommand{\Z}{\mathbb{Z}}

\newtheorem{thm1}{Theorem}
\newtheorem{thm2}[thm1]{Theorem}
\newtheorem{cor1}{Corollary}
\newtheorem{cor2}[cor1]{Corollary}
\newtheorem{smooth1}{Smooth Numbers Result}
\newtheorem{smooth2}[smooth1]{Smooth Numbers Result}
\newtheorem{smooth3}[smooth1]{Smooth Numbers Result}
\newtheorem{major1}{Major Arc Estimate}
\newtheorem{major2}[major1]{Major Arc Estimate}

\newtheorem{harmonic1}{Harmonic Analysis Result}
\newtheorem{harmonic2}[harmonic1]{Harmonic Analysis Result}
\newtheorem{prop1}{Proposition}
\newtheorem{prop2}[prop1]{Proposition}
\newtheorem{prop3}[prop1]{Proposition}
\newtheorem{prop4}[prop1]{Proposition}

\newtheorem{prop5}[prop1]{Proposition}

\newtheorem{trans1}{Transference Principle}

\begin{document}

\maketitle

\begin{abstract}
We investigate exponential sums over those numbers $\leq x$ all of whose prime factors are $\leq y$. We prove fairly good minor arc estimates, valid whenever $\log^{3}x \leq y \leq x^{1/3}$. Then we prove sharp upper bounds for the $p$-th moment of (possibly weighted) sums, for any real $p > 2$ and $\log^{C(p)}x \leq y \leq x$. Our proof develops an argument of Bourgain, showing this can succeed without strong major arc information, and roughly speaking it would give sharp moment bounds and restriction estimates for any set sufficiently factorable relative to its density.

By combining our bounds with major arc estimates of Drappeau, we obtain an asymptotic for the number of solutions of $a+b=c$ in $y$-smooth integers less than $x$, whenever $\log^{C}x \leq y \leq x$. Previously this was only known assuming the Generalised Riemann Hypothesis. Combining them with transference machinery of Green, we prove Roth's theorem for subsets of the $y$-smooth numbers, whenever $\log^{C}x \leq y \leq x$. This provides a deterministic set, of size $\approx x^{1-c}$, inside which Roth's theorem holds.
\end{abstract}

\section{Introduction}
One of the key tools in analytic number theory is the use of exponential sums over arithmetically interesting sets $\mathcal{S}$, namely sums of the form
$$ \sum_{\substack{n \leq x, \\ n \in \mathcal{S}}} e(n\theta) , $$
where $\theta \in \R$ and $e(z) := e^{2\pi i z}$ is the complex exponential. Such sums serve as Fourier transforms of $\mathcal{S}$, and are particularly useful in additive problems since one can detect solutions of additive equations very neatly in terms of exponential sums.

Exponential sums over primes are crucial in proving the ternary Goldbach conjecture, and also in proving results like Roth's theorem in the primes, which asserts that any positive density subset of the primes must contain a non-trivial three term arithmetic progression. Similarly, when $\mathcal{S}$ is the set of $k$-th powers such sums have been hugely studied in connection with Waring's problem. In this paper we investigate exponential sums over the set $\mathcal{S}(y)$ of $y$-smooth numbers, that is the set of numbers all of whose prime factors are less than or equal to $y$. Smooth numbers appear throughout analytic number theory, for example in proving some of the sharpest results on Waring's problem, and in analysing the performance of cryptographic algorithms. This is described in the survey articles of Granville~\cite{granville} and Hildebrand and Tenenbaum~\cite{ht3}.

\vspace{12pt}
Before giving a precise description of our results and their consequences, we try to indicate what is at stake by considering a particular additive problem, and the type of exponential sum information one would need to solve it.

Suppose we were seeking an asymptotic or lower bound for the number $N(x,y)$ of solutions to the equation $a+b=c$, with $a,b,c \leq x$ all $y$-smooth numbers. This is a much studied problem with connections to the $abc$ conjecture, as we will describe later. As usual in the circle method, the orthogonality of additive characters $e(n\theta)$ lets us write
$$ N(x,y) = \int_{0}^{1} \Biggl(\sum_{\substack{n \leq x, \\ n \in \mathcal{S}(y)}} e(n\theta) \Biggr)^{2} \overline{\Biggl(\sum_{\substack{n \leq x, \\ n \in \mathcal{S}(y)}} e(n\theta) \Biggr)} d\theta . $$
Let us also note that, heuristically, we expect to have $N(x,y) \asymp \Psi(x,y)^{3}/x$ on a very wide range of $y$, where $\Psi(x,y) := \sum_{\substack{n \leq x, \\ n \in \mathcal{S}(y)}} 1$. This is because there are $\Psi(x,y)^{3}$ unrestricted choices of $a,b,c$, and the condition that $a+b=c$ might be expected to hold with ``probability'' of order $1/x$ over all such choices (since $a,b,c$ have size at most $x$).

Now the usual way to proceed is to define a subset $\mathfrak{M} \subseteq [0,1]$ of {\em major arcs}, typically a union of small intervals around rational numbers $a/q$ with small denominators. One chooses $\mathfrak{M}$ such that one can obtain asymptotics for $\sum_{\substack{n \leq x, \\ n \in \mathcal{S}(y)}} e(n\theta)$ when $\theta \in \mathfrak{M}$. Moreover, one usually expects the largest values of the exponential sum to occur on the major arcs, and therefore that they will make the dominant contribution to the integral.

Let $\mathfrak{m} := [0,1] \backslash \mathfrak{M}$ denote the complementary set of {\em minor arcs}. In view of the above discussion, one would usually proceed by writing
\begin{eqnarray}
N(x,y) & = & \int_{\mathfrak{M}} \Biggl(\sum_{\substack{n \leq x, \\ n \in \mathcal{S}(y)}} e(n\theta) \Biggr)^{2} \overline{\Biggl(\sum_{\substack{n \leq x, \\ n \in \mathcal{S}(y)}} e(n\theta) \Biggr)} d\theta + O\Biggl(\sup_{\theta \in \mathfrak{m}} \Biggl|\sum_{\substack{n \leq x, \\ n \in \mathcal{S}(y)}} e(n\theta) \Biggr| \int_{0}^{1} \Biggl|\sum_{\substack{n \leq x, \\ n \in \mathcal{S}(y)}} e(n\theta) \Biggr|^{2} d\theta \Biggr) \nonumber \\
& = & \int_{\mathfrak{M}} \Biggl(\sum_{\substack{n \leq x, \\ n \in \mathcal{S}(y)}} e(n\theta) \Biggr)^{2} \overline{\Biggl(\sum_{\substack{n \leq x, \\ n \in \mathcal{S}(y)}} e(n\theta) \Biggr)} d\theta + O\Biggl(\sup_{\theta \in \mathfrak{m}} \Biggl|\sum_{\substack{n \leq x, \\ n \in \mathcal{S}(y)}} e(n\theta) \Biggr| \Psi(x,y) \Biggr) , \nonumber
\end{eqnarray}
where the second equality uses Parseval's identity to bound the mean square.

We now come to the crucial point. Parseval's identity, although simple and general, is also extremely powerful, since it shows that on average the exponential sum $\sum_{\substack{n \leq x, \\ n \in \mathcal{S}(y)}} e(n\theta)$ exhibits squareroot cancellation. Nevertheless, since we expect that $N(x,y)$ (and therefore also the integral over the major arcs) should have size $\asymp \Psi(x,y)^{3}/x$, in order to obtain a non-trivial result we would need to know (somewhat more than)
$$ \sup_{\theta \in \mathfrak{m}} \Biggl|\sum_{\substack{n \leq x, \\ n \in \mathcal{S}(y)}} e(n\theta) \Biggr| \ll \frac{\Psi(x,y)^{2}}{x} . $$
Here the trivial bound is $\Psi(x,y)$, so we ask for a saving of $\Psi(x,y)/x$. But if $y=\log^{K}x$, say, for any fixed $K \geq 1$, then it is known that $\Psi(x,\log^{K}x) = x^{1-1/K + o(1)}$ as $x \rightarrow \infty$, and so we need a {\em power saving} $x^{-1/K + o(1)}$ on the minor arcs. It is almost never possible to define a set of major arcs on which one can obtain asymptotics, and which is large enough that one has a power saving on the complementary minor arcs, except by assuming very strong conjectures like the Generalised Riemann Hypothesis. Note that if one were studying a ternary additive problem involving the primes, say, then one would only need to save a factor of about $1/\log x$ on the minor arcs, corresponding to the density of the primes. In this sense problems involving primes are much easier than problems involving very smooth numbers.

At this point, a circle method enthusiast might suggest various {\em pruning} manoeuvres to try to make progress. We will discuss this possibility later, but it has not yet made progress on the above problem, and our approach in this paper is different.

Note that, similarly as above, we can write
$$ N(x,y) = \int_{\mathfrak{M}} \Biggl(\sum_{\substack{n \leq x, \\ n \in \mathcal{S}(y)}} e(n\theta) \Biggr)^{2} \overline{\Biggl(\sum_{\substack{n \leq x, \\ n \in \mathcal{S}(y)}} e(n\theta) \Biggr)} d\theta + O\Biggl(\sup_{\theta \in \mathfrak{m}} \Biggl|\sum_{\substack{n \leq x, \\ n \in \mathcal{S}(y)}} e(n\theta) \Biggr|^{0.9} \int_{0}^{1} \Biggl|\sum_{\substack{n \leq x, \\ n \in \mathcal{S}(y)}} e(n\theta) \Biggr|^{2.1} d\theta \Biggr) , $$
say. Heuristically one might now expect that
$$ \int_{0}^{1} \Biggl|\sum_{\substack{n \leq x, \\ n \in \mathcal{S}(y)}} e(n\theta) \Biggr|^{2.1} d\theta \ll \max\left\{\frac{\Psi(x,y)^{2.1}}{x}, \Psi(x,y)^{1.05}\right\} , $$
since the first term reflects a contribution from major arcs, e.g. points at distance $\leq 1/x$ from 0, and the second term reflects squareroot cancellation. Provided $y \geq \log^{K}x$ for sufficiently large fixed $K$, the term $\Psi(x,y)^{2.1}/x$ will be the dominant one. Now {\em if} we could show that indeed $\int_{0}^{1} \left|\sum_{\substack{n \leq x, \\ n \in \mathcal{S}(y)}} e(n\theta) \right|^{2.1} d\theta \ll \frac{\Psi(x,y)^{2.1}}{x}$, then in order for the minor arc contribution to be small we would only need to show that
$$ \sup_{\theta \in \mathfrak{m}} \Biggl|\sum_{\substack{n \leq x, \\ n \in \mathcal{S}(y)}} e(n\theta) \Biggr|^{0.9} = o(\Psi(x,y)^{0.9}) . $$
In other words, we would no longer need to obtain a power saving, but rather {\em any} saving at all would be sufficient.

In this paper, we will prove sharp mean value estimates $\int_{0}^{1} \left|\sum_{\substack{n \leq x, \\ n \in \mathcal{S}(y)}} e(n\theta) \right|^{p} d\theta \ll \frac{\Psi(x,y)^{p}}{x}$ when $p > 2$. In fact, we will prove a much more general {\em restriction} result for arbitrarily weighted mean values $\int_{0}^{1} \left|\sum_{\substack{n \leq x, \\ n \in \mathcal{S}(y)}} a_{n} e(n\theta) \right|^{p} d\theta$. We will also prove some new pointwise minor arc bounds. These results have applications to ternary additive problems involving smooth numbers, and also to certain additive combinatorics questions.

We should emphasise that there is nothing new about trying to bound fractional moments of exponential sums. However, existing techniques would generally require an equivalent amount of major arc information as would make the standard circle method approach to ternary additive problems work directly. In our problems we do not have such information available. In contrast, we will obtain sharp mean value bounds in a quite general way, requiring very little major arc information.

\subsection{Brief survey of previous results}
Before explaining our new results, we indicate more precisely what was known about exponential sums over smooth numbers, with some comparisons to what is known about sums over primes.

\vspace{12pt}
{\em Major arc estimates.} If $\theta = a/q + \delta$, for some small $q \in \N$, some $(a,q)=1$, and some small real $\delta$, then $\theta$ is often described as lying on a {\em major arc}, and one can evaluate $\sum_{\substack{n \leq x, \\ n \in \mathcal{S}}} e(n\theta)$ quite precisely by splitting into arithmetic progressions to the small modulus $q$. In particular, if $\mathcal{S}$ is ``well distributed'' in arithmetic progressions to modulus $q$ then one can hope to obtain an asymptotic formula for $\sum_{\substack{n \leq x, \\ n \in \mathcal{S}}} e(n\theta)$.

For primes, one can show (see chapter 26 of Davenport~\cite{davenport}) that for any $A > 0$,
$$ \sum_{n \leq x} \Lambda(n) e(n(a/q+\delta)) = \frac{\mu(q)}{\phi(q)} T(\delta) + O_{A}((1+|\delta x|)x e^{-c\sqrt{\log x}}) \;\;\;\;\; \forall q \leq \log^{A}x, \; (a,q)=1 , $$
where $\mu(q)$ denotes the M\"{o}bius function, $\phi(q)$ denotes Euler's totient function, and $T(\delta) := \sum_{n \leq x} e(n\delta)$ is the sum of a geometric progression.

In the case of smooth numbers, Drappeau~\cite{drappeausommes} recently showed that if $\log^{c_{1}}x \leq y \leq x$, and $q \leq y^{c_{2}}$, and $(a,q)=1$, for certain absolute constants $c_{1},c_{2} > 0$, then
$$ \sum_{\substack{n \leq x, \\ n \in \mathcal{S}(y)}} e(n(a/q+\delta)) = \sum_{\substack{n \leq x, \\ n \in \mathcal{S}(y)}} \frac{\mu(q/(q,n))}{\phi(q/(q,n))} e(n\delta) + \text{error terms} , $$
and he gave useful estimates for the sum on the right hand side. We state these results carefully in $\S 2.2$, since we shall need to use them. Drappeau's estimates extend the permitted range of $y$ (and $q$ and $\delta$) as compared with previous results, for example those of La Bret\`{e}che~\cite{dlbLMS}, La Bret\`{e}che and Granville~\cite{dlbgranville}, and Fouvry and Tenenbaum~\cite{fouvryten0}.

\vspace{12pt}
{\em Minor arc estimates.} If $\theta$ isn't close enough to a rational with small denominator to be treated by major arc methods, then $\theta$ is usually described as lying on the {\em minor arcs}, and we usually settle for obtaining non-trivial upper bounds for $\sum_{\substack{n \leq x, \\ n \in \mathcal{S}}} e(n\theta)$.

For primes, the most common minor arc estimates (see chapter 25 of Davenport~\cite{davenport}) are of the following rough form, which is non-trivial provided $\log^{8}x \ll q \ll x/\log^{8}x$: 
$$ \left|\sum_{n \leq x} \Lambda(n) e(n(a/q+\delta)) \right| \ll x \log^{4}x \left(\frac{1}{\sqrt{q}} + \frac{1}{x^{1/5}} + \sqrt{\frac{q}{x}}\right) \;\;\;\;\; \forall q \in \N, \; (a,q)=1, |\delta| \leq \frac{1}{q^{2}} . $$
One can also prove stronger results with an explicit dependence on $\delta$ (e.g. Helfgott~\cite{helfgott} exploits such results in his work on the ternary Goldbach conjecture, building on their use elsewhere in the circle method). All these ultimately depend on formatting $\sum_{n \leq x} \Lambda(n) e(n(a/q+\delta))$ into multiple sums, as pioneered by Vinogradov, and then using the Cauchy--Schwarz inequality and completing one sum. The bound decays like $1/\sqrt{q}$ (and, in the $\delta$-dependent estimates, roughly like $1/\sqrt{q(1+|\delta x|)}$), but we expect decay roughly like $1/q(1+|\delta x|)$ on a wide range, as in the major arc estimates. This loss is ultimately due to the use of Cauchy--Schwarz, a point we shall return to later.

For smooth numbers, the best results that are comparable to the prime number estimates seem to be the following bounds, in which $q \in \N$ and $(a,q)=1$:
$$ \left| \sum_{\substack{n \leq x, \\ n \in \mathcal{S}(y)}} e(n(\frac{a}{q}+\delta)) \right| \ll \left\{ \begin{array}{ll}
      x (1+|\delta x|) \log^{3}x \left(\frac{\sqrt{y}}{x^{1/4}} + \frac{1}{\sqrt{q}} + \sqrt{\frac{qy}{x}} \right)& \text{if} \; 3 \leq y \leq \sqrt{x}, \\
     x\sqrt{\log x} \log y \left(\sqrt{\frac{y}{x}} + \frac{1}{\sqrt{q}} + \sqrt{\frac{q}{x}} + e^{-\sqrt{\log x}} \right)  & \text{if} \; x,y \geq 2, \; |\delta| \leq \frac{1}{q^{2}} .
\end{array} \right. $$
These are due to Fouvry and Tenenbaum~\cite{fouvryten0} and to La Bret\`{e}che~\cite{dlbLMS}, respectively. They are a bit unsatisfactory unless $y$ and $q$ are fairly large, because the leading term is $x$ rather than the trivial bound of $\Psi(x,y) := \sum_{\substack{n \leq x, \\ n \in \mathcal{S}(y)}} 1$, which may be {\em much} smaller.

Our first theorem, stated shortly, will be a minor arc estimate where the leading term is $\Psi(x,y)$, valid on the wide range $\log^{3}x \leq y \leq x^{1/3}$, and with an explicit dependence on $\delta$ that further strengthens the bound on certain ranges.

\vspace{12pt}
{\em Mean value theorems.} For many sets $\mathcal{S}$ we expect a mean value bound
$$ \int_{0}^{1} \Biggl|\sum_{\substack{n \leq x, \\ n \in \mathcal{S}}} e(n\theta) \Biggr|^{p} d\theta \ll_{p} \max\Biggl\{\frac{1}{x} \Biggl(\sum_{\substack{n \leq x, \\ n \in \mathcal{S}}} 1 \Biggr)^{p} , \Biggl(\sum_{\substack{n \leq x, \\ n \in \mathcal{S}}} 1 \Biggr)^{p/2}\Biggr\} $$
for any real $p \geq 2$, since the first term reflects a contribution from major arcs, e.g. points at distance $\leq 1/x$ from 0, and the second term reflects squareroot cancellation. Note that when $p=2$ the second term is always the dominant one, but for larger $p$ the first term may dominate.

For primes, it is fairly easy to show that indeed $\int_{0}^{1} \left|\sum_{n \leq x} \Lambda(n) e(n\theta) \right|^{p}d\theta \ll_{p} x^{p-1}$ whenever $p > 2$, simply by adapting the argument we outlined when thinking about $N(x,y)$. For if $\mathfrak{M}$ denotes a suitable set of major arcs (whose definition will depend on $p$), and $\mathfrak{m} := [0,1]\backslash \mathfrak{M}$ the complementary set of minor arcs, then
$$ \int_{0}^{1} \left|\sum_{n \leq x} \Lambda(n) e(n\theta) \right|^{p}d\theta \ll \int_{\mathfrak{M}} \left|\sum_{n \leq x} \Lambda(n) e(n\theta) \right|^{p} + \sup_{\theta \in \mathfrak{m}} \left|\sum_{n \leq x} \Lambda(n) e(n\theta) \right|^{p-2} x \log x , $$
by Parseval's identity. One can evaluate the integral over $\mathfrak{M}$ using pointwise asymptotics, and one has bounds for the supremum on the minor arcs that save at least a factor $1/\log^{1/(p-2)}x$. The crucial point here is that the primes are logarithmically far from being a dense set, producing the unwanted factor $\log x$, and one has major arc estimates on a wide enough range of $q$ and $\delta$ to compensate for that logarithmic loss.

For smooth numbers, when $C \sqrt{\log x} \log\log x \leq \log y$ one could prove a sharp bound $\int_{0}^{1} \left|\sum_{\substack{n \leq x, \\ n \in \mathcal{S}(y)}} e(n\theta)\right|^{p} d\theta \ll_{p} \Psi(x,y)^{p}/x$ for any $p > 2$, as for the primes. (See $\S 2.1$ for some explanation of this range.) But for smaller $y$ the known major arc estimates are valid on too small a range of $q$ to compensate for the very low density of the $y$-smooth numbers, even when accompanied by our new minor arc estimates. Thus it has been impossible to handle many additive problems involving smooth numbers when $y$ is small.

Our second and main theorem will be a sharp upper bound for $\int_{0}^{1} \left|\sum_{\substack{n \leq x, \\ n \in \mathcal{S}(y)}} e(n\theta)\right|^{p} d\theta$, that is valid whenever $p > 2$ and $\log^{C(p)}x \leq y \leq x$.

\vspace{12pt}
{\em Restriction theory.} Lastly we reach a class of results that may be less familiar to some analytic number theorists, namely {\em majorant} and {\em restriction} theorems. This terminology is to some extent imported from Euclidean harmonic analysis, but all we mean is that we would like to show, for any complex numbers $|a_{n}| \leq 1$ and any $p > 2$, and for a given nice set $\mathcal{S}$, that
$$ \int_{0}^{1} \left|\sum_{\substack{n \leq x, \\ n \in \mathcal{S}}} a_{n} e(n\theta) \right|^{p}d\theta \ll_{p} \int_{0}^{1} \left|\sum_{\substack{n \leq x, \\ n \in \mathcal{S}}} e(n\theta) \right|^{p}d\theta , \;\;\; \text{or} \;\;\; \int_{0}^{1} \left|\sum_{\substack{n \leq x, \\ n \in \mathcal{S}}} a_{n} e(n\theta) \right|^{p}d\theta \ll_{p} \left(\sum_{\substack{n \leq x, \\ n \in \mathcal{S}}} 1 \right)^{p}/x . $$
This is another mean value problem, but the point is that we would like results that don't require the $a_{n}$ to be ``nice'', but hold simply because $\mathcal{S}$ has good properties.

For primes, the restriction estimate
$$ \int_{0}^{1} \left| \sum_{p \leq x} a_{p} e(p\theta) \right|^{q} d\theta \ll_{q} \pi(x)^{q}/x \;\;\;\;\; \forall q > 2 $$
was first proved by Bourgain~\cite{bourgain} (though expressed rather differently, and with an $L^{2}$ rather than an $L^{\infty}$ condition on the coefficients $a_{p}$), and used to show that random subsets of the primes have certain nice harmonic analysis properties. We will discuss his proof later. These results gained great prominence when Green~\cite{greenroth} gave another proof, and used the restriction estimate to prove Roth's theorem in the primes.

For smooth numbers, the author is not aware of any previous restriction or majorant results, though they could have been proved for large $y$ (e.g. for $\log y \geq C \sqrt{\log x} \log\log x$) by direct use of Bourgain's~\cite{bourgain} method. Green's~\cite{greenroth} method crucially exploits the fact that primes can be majorised efficiently by sieve processes\footnote{The characteristic function of the primes is decomposed into pieces such that, roughly speaking, earlier pieces are dense (so Parseval-type arguments are efficient) but not particularly Fourier uniform, whereas later pieces are sparser (like the logarithmic sparseness of the primes) but only produce small exponential sums. This crucial decomposition comes from sieving with increasing levels.}, and since that is well known to fail for smooth numbers his method doesn't seem applicable.

We shall prove a restriction theorem for the $y$-smooth numbers whenever $p > 2$ and $\log^{C(p)}x \leq y \leq x$. The key new point is that we can do this even though pointwise major and minor arc estimates are insufficient even to sharply bound the {\em unweighted} mean value. Indeed, we deduce our mean value theorem from our restriction theorem.

\subsection{Statement of new results and applications}
Firstly we shall prove a minor arc estimate that is sensitive to the density of the smooth numbers. This result involves a standard quantity $\alpha(x,y)$, whose definition and properties will be summarised in $\S 2.1$. For the moment, we only say very roughly that $\alpha(x,y)$ is ``close to'' 1 provided $(\log y)/\log\log x$ is large.
\begin{thm1}
Let $\log^{3}x \leq y \leq x^{1/3}$ be large. Suppose that $\theta = \frac{a}{q} + \delta$ for some $(a,q)=1$ and some $\delta \in \R$. Provided that $q^{2}y^{3}(1+|\delta x|)^{2} \leq x/4$, we have the bound
$$ \Biggl| \sum_{\substack{n \leq x, \\ n \in \mathcal{S}(y)}} e(n\theta) \Biggr| \ll \frac{\Psi(x,y)}{\sqrt{q(1+|\delta x|)}} (q(1+|\delta x|))^{(3/2)(1-\alpha(x,y))} u^{3/2} \log u \log x \sqrt{\log(2+|\delta x|) \log(qy)} , $$
where $\Psi(x,y) := \sum_{\substack{n \leq x, \\ n \in \mathcal{S}(y)}} 1$, where $u:= (\log x)/\log y$, and where $\alpha(x,y)$ denotes the saddle-point corresponding to the $y$-smooth numbers less than $x$.
\end{thm1}

For a reader who is less familiar with smooth numbers, we emphasise again that when $y$ is small they are very sparse, for example $\Psi(x,\log^{K}x)=x^{1-1/K + o(1)}$ for any fixed $K \geq 1$, as $x \rightarrow \infty$. Thus it is important to have estimates in which the leading term is $\Psi(x,y)$ rather than $x$. The condition that $q^{2}y^{3}(1+|\delta x|)^{2} \leq x/4$, which arises naturally in the proof, is not really restrictive, since we can generally assume that $y \leq x^{1/100}$ and (after applying Dirichlet's approximation theorem to $\theta$) that $q \leq x^{0.55}$ and $|\delta x| \leq x^{0.45}/q$, say. Then the only case we cannot handle with Theorem 1 is where $x^{0.48} \leq q \leq x^{0.55}$, for which we can use Fouvry and Tenenbaum's~\cite{fouvryten0} minor arc estimate (stated earlier) and obtain a substantial saving because $q$ is so large. Note also that, because our estimate improves as $\delta$ increases, we don't need to restrict to narrow arcs around $a/q$ even when $q$ is small, which is helpful in applications.

The proof of Theorem 1 isn't very difficult, since it works in the standard way of formatting the exponential sum into double sums, which is easy because the smooth numbers have good factorisation properties, and then using the Cauchy--Schwarz inequality. The only issue that arises is obtaining a reasonable upper bound for the quantity of $y$-smooth numbers in a segment of an arithmetic progression, on a wide range of $y$. By combining a method of Friedlander~\cite{friedlanderupper} with modern smooth number estimates we obtain upper bounds that aren't too bad (see $\S 2.1$, below), and by controlling the lengths of our double sums we can ensure the losses are acceptable.

Our second and main result is the following restriction theorem. Note that by choosing all the coefficients $a_{n}$ to be 1 we obtain a mean value result for exponential sums over smooth numbers, which was itself unknown for small $y$.
\begin{thm2}
There exists an absolute constant $C > 0$ such that the following is true.

Let $p > 2$, suppose $x$ is large enough in terms of $p$, and let $\log^{C\max\{1,1/(p-2)\}}x \leq y \leq x$. Then for any complex numbers $(a_{n})_{n \leq x}$ with absolute values at most 1, we have
$$ \int_{0}^{1} \Biggl| \sum_{\substack{n \leq x, \\ n \in \mathcal{S}(y)}} a_{n} e(n\theta) \Biggr|^{p} d\theta \ll_{p} \frac{\Psi(x,y)^{p}}{x} . $$
\end{thm2}

The lower bound condition $y \geq \log^{C\max\{1,1/(p-2)\}}x$ is sharp up to the value of $C$, at least for $p$ close to 2. For it is easy to show that $\int_{0}^{1} \left| \sum_{\substack{n \leq x, \\ n \in \mathcal{S}(y)}} e(n\theta) \right|^{p} d\theta \gg_{p} \Psi(x,y)^{p/2}$ whenever $p \geq 2$, and since $\Psi(x,\log^{C/(p-2)}x) = x^{1-(p-2)/C+o(1)}$ we are close to the transition point where $\Psi(x,y)^{p/2} > \Psi(x,y)^{p}/x$.

To explain the proof of Theorem 2, we first explain Bourgain's~\cite{bourgain} restriction argument for the primes. Bourgain~\cite{bourgain} bounds the measure of the set of $\theta$ for which $\left| \sum_{p \leq x} a_{p} e(p\theta) \right| \geq \delta \pi(x)$, which will yield the restriction theorem if one can obtain good bounds for a sufficient range of $\delta$. He makes a clever application of the Cauchy--Schwarz inequality, (described as ``linearization'', and which the reader might recognise from Hal\'{a}sz--Montgomery type arguments for Dirichlet polynomials), after which the unknown coefficients $a_{p}$ are removed and one can insert bounds for {\em unweighted} exponential sums over primes, with $\theta$ replaced by differences $\theta_{r}-\theta_{s}$. Since one has good major arc estimates for primes, one obtains bounds for the set of $\theta$ after finally exploiting the good spacing of the major arcs (as encoded in a result like Harmonic Analysis Result 2, in $\S 2.3$ below).

We will show that Bourgain's argument can be run without inserting strong major arc estimates. In fact, we will show that his Cauchy--Schwarz ``duality'' step combines very naturally with the usual {\em proof} of minor arc estimates, that involves producing double sums. The key point is that, when proving pointwise bounds like Theorem 1, the use of the Cauchy--Schwarz inequality is essential but wasteful, losing a squareroot factor. In contrast, in Bourgain's argument we already have a Cauchy--Schwarz step that is not wasteful, since it involves sets of $\theta$ rather than $\mathcal{S}$. It turns out that this one Cauchy--Schwarz step can suffice for everything, so that the double sums approach to minor arcs, which loses a squareroot factor when seeking pointwise bounds, becomes quite sharp when seeking mean value or restriction results. This approach is not at all restricted to smooth numbers, and (though we don't attempt a general formulation) a rough indication of the requirement on $\mathcal{S} \subseteq [1,x]$ is the following: for any $1 \leq K \leq (\frac{x}{\#\mathcal{S}} \log x)^{100}$, one must be able to find some set $\mathcal{S}_{K} \subseteq [1,x/K]$ such that
$$ \mathcal{S} \subseteq \{mn : m \in \mathcal{S}_{K}, n \leq x/m\}, \;\;\;\;\; \text{and} \;\;\;\;\; \sum_{m \in \mathcal{S}_{K}} \sum_{n \leq x/m} 1 \ll (\#\mathcal{S}) K^{0.01} \log^{1000}x . $$
(Here the numbers 100, 0.01, 1000, are indicative only, and will really depend on how close $p$ is to 2 and on how much major arc information one has about $\mathcal{S}$.)

We actually use Theorem 1, combined with the Erd\H{o}s--Tur\'{a}n inequality, in the proof of Theorem 2, to reduce to a situation where we only need look at differences $\theta_{r}-\theta_{s}$ that are ``fairly close'' to a rational $a/q$ with ``fairly small'' denominator (a reduction that lets us set the lengths of double sums properly). We also use sharp major arc estimates, but only in the very restricted case where $q$ has size $\log^{O(1)}x$, because our main argument loses a few logarithmic factors and we need a different approach if the other terms aren't large enough to compensate for this.

Let us mention a few other restriction results in the number theoretic literature. Green and Tao~\cite{greentaorest} prove a restriction theorem for certain sets that are well controlled by the sieve, by a neat argument that is something of a hybrid between the approaches of Bourgain~\cite{bourgain} and Green~\cite{greenroth}. Recently, Keil~\cite{keil1,keil2} has used another hybrid approach to study solutions of quadratic systems and forms in dense variables. He decomposes the counting function corresponding to his quadratic system, but does this based on a major arc decomposition on the Fourier side, rather than with sieves. This again appears to require sharp major arc estimates on a wide range. We also note the techniques of {\em pruning}, which are used to great effect in the circle method when limited major arc information is available. In fact, results like Br\"{u}dern's pruning lemma (see Lemma 2 of \cite{brudern}) seem notably close to results from restriction theory, such as the technique of Bourgain~\cite{bourgain} that we quote as Harmonic Analysis Result 2, below. But to apply pruning one generally needs at least some exponential sums around for which good major arc information is available, and in our problems we have none.

\vspace{12pt}
Finally turning to applications, by combining Theorems 1 and 2 with major arc estimates of Drappeau~\cite{drappeausommes} we can readily handle ternary additive problems involving smooth numbers. For example, we obtain the following result.
\begin{cor1}
There exists a large absolute constant $K > 0$ such that, for any large $\log^{K}x \leq y \leq x$, (and writing $u:= (\log x)/\log y$),
$$ \#\{(a,b,c) \in \mathcal{S}(y)^{3} : a,b,c \leq x, \; a+b=c\} = \frac{\Psi(x,y)^{3}}{2x} \left(1+O\left(\frac{\log(u+1)}{\log y}\right) \right) . $$
\end{cor1}
The error term $\log(u+1)/\log y$ tends to zero if $(\log y)/\log\log x \rightarrow \infty$, and for smaller $y$ (still satisfying $y \geq \log^{K}x$) one can obtain an error term tending to zero by replacing the main term $\Psi(x,y)^{3}/2x$ by something a little more complicated, but still explicit.

Assuming the Generalised Riemann Hypothesis, Lagarias and Soundararajan~\cite{lagariassound0, lagariassound} (see also \cite{drappeauabc}) have previously proved Corollary 1 for all $\log^{8+\epsilon}x \leq y \leq e^{\log^{1/2-\epsilon}x}$. There is analogous work of Ha~\cite{ha} in the function field setting. Lagarias and Soundararajan investigated this equation as an analogue of the $abc$ conjecture, with the smoothness bound $y$ taking the place of the radical of $abc$ in that conjecture. Corollary 1 now proves unconditionally that Lagarias and Soundararajan's ``$xyz$-smoothness exponent'' is at most $K$, and in particular is finite\footnote{Strictly speaking, Lagarias and Soundararajan need $a,b,c$ to be coprime, but a simple inclusion-exclusion argument as in $\S 8$ of their paper~\cite{lagariassound} can be used to impose that condition.}.

Previously, the best known unconditional result like Corollary 1 was due to Drappeau~\cite{drappeausommes}, who could handle the range $e^{c\sqrt{\log x}\log\log x} \leq y \leq x$, improving on La Bret\`{e}che and Granville~\cite{dlbgranville} who could handle $e^{\log^{2/3+\epsilon}x} \leq y \leq x$. We also mention an older result of Balog and A. S\'{a}rk\H{o}zy~\cite{balogsarkozy}, who showed the existence of solutions to the equation $a+b+c=N$ with $a,b,c$ all $e^{3\sqrt{\log N \log\log N}}$-smooth, for any large $N$. They used exponential sums and the circle method, but with weights that allowed stronger estimates at the cost of losing an asymptotic for the number of solutions. G. S\'{a}rk\H{o}zy~\cite{gsarkozy} improved this a bit, replacing the constant 3 in the exponent by $\sqrt{3/2} + \epsilon$, using the large sieve.

\vspace{12pt}
Our second application is to additive combinatorics, where we combine Theorems 1 and 2, major arc estimates, and the transference machinery of Green~\cite{greenroth} to prove a version of Roth's theorem for smooth numbers.
\begin{cor2}
For any $\beta > 0$ there exists a large constant $K = K(\beta) > 0$ such that the following is true. For any $x$ that is large enough in terms of $\beta$, any $\log^{K}x \leq y \leq x$, and any set $B \subseteq \mathcal{S}(y) \cap [1,x]$ such that $\#B \geq \beta \Psi(x,y)$, there exist integers $b$ and $d \neq 0$ such that
$$ b, b+d, b+2d \in B . $$
\end{cor2}

As we mentioned, $\Psi(x,\log^{K}x)=x^{1-1/K+o(1)}$ as $x \rightarrow \infty$, so Corollary 2 exhibits a very sparse set inside which Roth's theorem holds. So far as the author knows, the only comparably sparse deterministic sets inside which Roth's theorem is known come from recent work of Mirek~\cite{mirek}, who handled sets like the Piatetski--Shapiro primes (that is primes of the form $\lfloor n^{1/\gamma} \rfloor$, where $\gamma$ is sufficiently close to 1 and $\lfloor \cdot \rfloor$ denotes integer part). In those cases the form of the function $n^{1/\gamma}$ again lets one obtain strong pointwise exponential sum bounds, unlike in our problem. Despite this paucity of explicit examples, Kohayakawa, {\L}uczak and R\"{o}dl~\cite{klr} proved that, with probability tending to 1 as $x \rightarrow \infty$, Roth's theorem holds inside random subsets of $[1,x]$ of cardinality $x^{1/2+o(1)}$ (roughly speaking).

To prove a result like Corollary 2, one needs a sharp restriction estimate for some $2 < p < 3$ for exponential sums over the set of interest (in this case the smooth numbers), and one needs to find a function $\nu(n)$ that upper bounds the characteristic function of that set, whose total sum is within a constant factor of the size of the set, {\em and} such that all exponential sums $\sum \nu(n) e(an/N)$ at certain well spaced points $a/N$ are small. We give a precise statement of what is required as Transference Principle 1, in $\S 6$. In our case, Theorem 2 immediately provides the restriction estimate, and the obvious choice of $\nu(n)$ as the characteristic function $\textbf{1}_{n \in \mathcal{S}(y)}$ is acceptable.

The reader might object that Corollary 2 doesn't really prove Roth's theorem inside the $\log^{K}x$-smooth numbers for any fixed $K$, because $K$ is required to grow with $\beta$. One could prove Corollary 2 with a large fixed $K$, but the proof would become more complicated because one would need to compensate for some irregularities of distribution of very smooth numbers, both $p$-adically (they are disproportionately unlikely to be coprime to small primes) and in an Archimedean sense (there are more $\log^{K}x$-smooth numbers between 1 and $x/2$ than between $x/2$ and $x$, say). One could handle this using an analogue of Green's~\cite{greenroth} {\em $W$-trick}, in which he restricted to working inside certain arithmetic progressions to overcome the irregular distribution of primes to small moduli (they are disproportionately likely to be coprime to the modulus). However, we have chosen not to present such an argument to avoid making this paper even longer.

\vspace{12pt}
We end this very long introduction by explaining the organisation of the rest of the paper. In $\S 2$ we state a few background results we shall need, concerning the distribution of smooth numbers, major arc estimates for smooth numbers, and a couple of harmonic analysis facts. Most of these are imported from the literature, but we need to prove a satisfactory upper bound for smooth numbers in segments of arithmetic progressions, and we also need to extend some known major arc estimates a little. Proofs of the new parts of major arc estimates are deferred to the appendix. In $\S 3$ we prove Theorem 1, and in $\S 4$ we prove Theorem 2. Finally, in $\S\S 5-6$ we deduce our two corollaries.

\section{Some tools}

\subsection{General results on smooth numbers}
In this subsection we state some general smooth number estimates, but not exponential sum estimates, that we shall need.

We begin with a celebrated result of Hildebrand and Tenenbaum~\cite{ht}, that gives an asymptotic for $\Psi(x,y) := \sum_{n \leq x, n \in \mathcal{S}(y)} 1$ on a very wide range of $x,y$, in terms of a {\em saddle point} $\alpha = \alpha(x,y)$. We won't actually need this result much, but the subsequent result (which we will use extensively) will make much more sense if the reader is familiar with Hildebrand and Tenenbaum's asymptotic.
\begin{smooth1}[Hildebrand and Tenenbaum, 1986]
We have uniformly for $x \geq y \geq 2$,
$$ \Psi(x,y) = \frac{x^{\alpha} \zeta(\alpha,y)}{\alpha \sqrt{2\pi(1+(\log x)/y)\log x \log y}} \left(1 + O\left(\frac{1}{\log(u+1)} + \frac{1}{\log y} \right) \right), $$
where $u = (\log x)/\log y$, $\zeta(s,y) := \sum_{n: n \textrm{ is } y \textrm{ smooth}} 1/n^{s} = \prod_{p \leq y}(1-p^{-s})^{-1}$ for $\Re(s) > 0$, and $\alpha = \alpha(x,y) > 0$ is defined by
$$ \sum_{p \leq y} \frac{\log p}{p^{\alpha}-1} = \log x. $$
\end{smooth1}

Hildebrand and Tenenbaum~\cite{ht} also established a simple approximation for $\alpha(x,y)$ on the whole range $2 \leq y \leq x$. Their Lemma 2 implies, in particular, that when $\log x < y \leq x$ one has
\begin{equation}\label{alphaapprox}
\alpha(x,y) = 1 - \frac{\log(u\log(u+1))}{\log y} + O\left(\frac{1}{\log y}\right).
\end{equation}
We also remark that, by definition, $\alpha(x,y)$ is a decreasing function of $x$ for any fixed $y$.

To orient an unfamiliar reader, we also provide a more explicit estimate for $\Psi(x,y)$ on a slightly smaller range. Hildebrand~\cite{hildebrand} showed that
$$ \Psi(x,y) = x\rho(u)\left(1+O_{\epsilon}\left(\frac{\log(u+1)}{\log y}\right) \right), \;\;\;\;\; e^{(\log\log x)^{5/3+\epsilon}} \leq y \leq x, $$
where the Dickman function $\rho(u)$ is a certain continuous function that satisfies $\rho(u) = e^{-(1+o(1))u\log u}$ as $u \rightarrow \infty$. Thus the $y$-smooth numbers are very sparse when $u=(\log x)/\log y$ is large\footnote{When $\log^{1+\epsilon}x \leq y \leq x$, Theorem 2(ii) and Corollary 2 of Hildebrand and Tenenbaum~\cite{ht} also imply that $\Psi(x,y) = x\rho(u)e^{O_{\epsilon}(u)}$. This implies the estimate $\Psi(x,\log^{K}x)=x^{1-1/K+o(1)}$ for any fixed $K > 1$, as $x \rightarrow \infty$, which we mentioned several times in the introduction (and which can actually be proved more simply, and is true when $K=1$ as well).}. In particular, when $y$ decreases below $e^{\sqrt{(1/2) \log x \log\log x}}$ there is a change in that $\Psi(x,y)/x \approx \rho(u)$ becomes smaller than $1/y$. As major arc estimates for smooth numbers are valid when $q \leq y^{c}$, this explains why a classical approach to ternary additive problems with smooth numbers fails for $y$ smaller than about $e^{\sqrt{(1/2) \log x \log\log x}}$.

\vspace{12pt}
The next result will be very important to our arguments. It is a neat ``local'' result of La Bret\`{e}che and Tenenbaum~\cite{dlbten}, that compares $\Psi(x,y)$ with $\Psi(x/d,y)$ on a very wide range of the parameters.
\begin{smooth2}[See Th\'{e}or\`{e}me 2.4(i) of La Bret\`{e}che and Tenenbaum~\cite{dlbten}]
Let $2 \leq y \leq x$, and suppose that $d \geq 1$. Then we have
$$ \Psi(x/d,y) \ll \frac{1}{d^{\alpha}} \Psi(x,y) , $$
where $\alpha = \alpha(x,y)$ denotes the saddle-point corresponding to the $y$-smooth numbers less than $x$.
\end{smooth2}

As mentioned in the introduction, we will also need an upper bound for the quantity of smooth numbers in a segment of an arithmetic progression. By combining Smooth Numbers Result 2 with a simple but powerful method of Friedlander~\cite{friedlanderupper}, we will prove the following result.
\begin{smooth3}
Let $\log X \leq y \leq X$ be large, and suppose $q \geq 1$ and $qy \leq Z \leq X$. Then
$$ \sum_{\substack{X \leq n \leq X+Z, \\ n \equiv a \; \text{mod} \; q, \\ n \in \mathcal{S}(y)}} 1 \ll \frac{Z}{qX} \Psi(X,y) \left(\frac{Xq}{Z}\right)^{1-\alpha} \log X , $$
where $\alpha = \alpha(X,y)$ denotes the saddle-point corresponding to the $y$-smooth numbers less than $X$.
\end{smooth3}

Note that if the $y$-smooth numbers were roughly equidistributed in short intervals and arithmetic progressions, one might expect an upper bound $\ll (Z/qX) \Psi(X,y)$, and if $y$ is at least a large power of $\log X$ then $1-\alpha$ will be ``close'' to 0, so the bound in Smooth Numbers Result 3 will be of roughly the correct shape (although a more precise bound would certainly be desirable). One could reduce or remove the factor $\log X$, and weaken the assumptions on $y$ and $Z$, with more work and by slightly reformulating the conclusion, but this will be unnecessary for us. Let us also remark that Friedlander~\cite{friedlanderupper} proved an upper bound roughly comparable to Smooth Numbers Result 3 for $y \geq e^{\log^{4/5}X}$ (and for $Z=X$), and later Balog and Pomerance~\cite{balogpomerance} used his method to handle $y \geq \exp\{(\log\log X)^{2}\}$, although their bound was weaker (primarily since they did not have access to a suitable result like Smooth Numbers Result 2).

\vspace{12pt}
Turning to the proof of Smooth Numbers Result 3, suppose at first (for simplicity) that $(a,q)=1$. Write $P(m)$ for the largest prime factor of the integer $m$, and write $p(n)$ for the smallest prime factor of the integer $n$. Then we see
$$ \sum_{\substack{X \leq n \leq X+Z, \\ n \equiv a \; \text{mod} \; q, \\ n \in \mathcal{S}(y)}} 1 = \sum_{\substack{Z/qy < m \leq Z/q, \\ m/P(m) \leq Z/qy, \\ m \in \mathcal{S}(y), \\ (m,q)=1}} \sum_{\substack{X/m \leq n \leq (X+Z)/m, \\ p(n) \geq P(m), \\ mn \equiv a \; \text{mod} \; q, \\ n \in \mathcal{S}(y)}} 1 , $$
since every smooth number from the interval $[X,X+Z]$ can be written uniquely as such a product $mn$, by revealing its prime factors one at a time, starting with the smallest, until the product exceeds $Z/qy$. Note that we may impose the condition that $(m,q)=1$ since otherwise the condition $mn \equiv a \; \text{mod} \; q$ cannot possibly be satisfied, given our assumption that $(a,q)=1$.

Next, if we write $m^{-1}$ to mean the multiplicative inverse of $m$ modulo $q$, and relax the inner summation, we obtain
$$ \sum_{\substack{X \leq n \leq X+Z, \\ n \equiv a \; \text{mod} \; q, \\ n \in \mathcal{S}(y)}} 1 \leq \sum_{\substack{Z/qy < m \leq Z/q, \\ m/P(m) \leq Z/qy, \\ m \in \mathcal{S}(y), \\ (m,q)=1}} \sum_{\substack{X/m \leq n \leq (X+Z)/m, \\ n \equiv a m^{-1} \; \text{mod} \; q}} 1 \leq \frac{2Z}{q} \sum_{\substack{Z/qy < m \leq Z/q, \\ m/P(m) \leq Z/qy, \\ m \in \mathcal{S}(y), \\ (m,q)=1}} \frac{1}{m} , $$
in view of the fact that $Z/m \geq q$. Dropping the restrictions that $m/P(m) \leq Z/qy$ and $(m,q)=1$, we see
\begin{eqnarray}
\sum_{\substack{Z/qy < m \leq Z/q, \\ m \in \mathcal{S}(y)}} \frac{1}{m} & \leq & \sum_{0 \leq j \leq (\log y)/\log 2} \frac{1}{2^{j}(Z/qy)} \sum_{\substack{2^{j}(Z/qy) \leq m \leq 2^{j+1}(Z/qy), \\ m \in \mathcal{S}(y)}} 1 \nonumber \\
& \leq & \sum_{0 \leq j \leq (\log y)/\log 2} \frac{1}{2^{j}(Z/qy)} \Psi\left(\frac{2^{j+1}Z}{Xqy} X, y \right) \nonumber \\
& \ll & \Psi(X,y) \left(\frac{Z}{Xqy}\right)^{\alpha} \sum_{0 \leq j \leq (\log y)/\log 2} \frac{2^{j(\alpha - 1)}}{(Z/qy)}, \nonumber
\end{eqnarray}
in view of Smooth Numbers Result 2. Finally, if we write $u = (\log X)/\log y$ then, since $\alpha = 1 - \frac{\log(u\log(u+1)) + O(1)}{\log y}$ when $\log X \leq y \leq X$, and since $y^{1-\alpha} = O(u\log(u+1))$, we certainly have
$$ \sum_{\substack{Z/qy < m \leq Z/q, \\ m \in \mathcal{S}(y)}} \frac{1}{m} \ll \frac{1}{X} \Psi(X,y) \left(\frac{Xqy}{Z}\right)^{1-\alpha} \frac{\log y}{\log(u+1)} \ll \frac{1}{X} \Psi(X,y) \left(\frac{Xq}{Z}\right)^{1-\alpha} \log X , $$
and Smooth Numbers Result 3 follows on putting everything together.

It only remains to remove our simplifying assumption that $(a,q)=1$. If $(a,q)=h$ for some $h$ that is not $y$-smooth then the sum in Smooth Numbers Result 3 has no terms, so the result is trivial. If $(a,q)=h$ for some $h \neq 1$ that is $y$-smooth then
$$ \sum_{\substack{X \leq n \leq X+Z, \\ n \equiv a \; \text{mod} \; q, \\ n \in \mathcal{S}(y)}} 1 = \sum_{\substack{X/h \leq n \leq (X+Z)/h, \\ n \equiv a/h \; \text{mod} \; q/h, \\ n \in \mathcal{S}(y)}} 1 , $$
where $(a/h,q/h)=1$, and where all the hypotheses of Smooth Numbers Result 3 are satisfied by the sum on the right. (Note that $y \leq X/q \leq X/h$, by assumption.) The result then follows by the above argument, and by using Smooth Numbers Result 2 again to obtain that $\Psi(X/h,y) \ll h^{-\alpha} \Psi(X,y)$.
\begin{flushright}
Q.E.D.
\end{flushright}

\subsection{Major arc estimates}
In this subsection we state some major arc estimates for exponential sums over smooth numbers. We will need these both as an ingredient for the proof of Theorem 2, to avoid logarithmic losses, and when deducing Corollaries 1 and 2.

\begin{major1}
There exist two absolute positive constants $c_{1},c_{2}$ such that the following is true. For any large $\log^{c_{1}}x \leq y \leq x$, any $q \leq y^{c_{2}}$, any $(a,q)=1$, any $\delta \in \R$, and any $A > 0$, we have
$$ \sum_{\substack{n \leq x, \\ n \in \mathcal{S}(y)}} e(n(a/q+\delta)) = V(x,y;q,\delta) + O_{A}(\Psi(x,y)(1+|\delta x|)( y^{-c_{2}} + \frac{e^{-c_{2}u/\log^{2}(u+1)}}{\log^{A}x})) , $$
where $u := (\log x)/\log y$, $V(x,y;q,\delta) := \sum_{\substack{n \leq x, \\ n \in \mathcal{S}(y)}} \frac{\mu(q/(q,n))}{\phi(q/(q,n))} e(n\delta)$, and $\mu, \phi$ denote the M\"{o}bius and Euler totient functions.
\end{major1}

Major Arc Estimate 1 reduces the study of exponential sums over smooth numbers on the major arcs to an investigation of the main term $V(x,y;q,\delta)$. The estimate more or less follows from Proposition 2.2 of Drappeau~\cite{drappeausommes}, with the choice $Q := \min\{y^{c_{2}},e^{c_{2}\sqrt{\log x}}\}$, except that proposition includes an extra term $W(x,y;q,\delta)$ reflecting the possible existence of an exceptional (Siegel) zero. Given our current knowledge about character sums over smooth numbers, it actually isn't too difficult to bound $W(x,y;q,\delta)$ satisfactorily on the complete range in Major Arc Estimate 1. We supply the relevant argument, which combines with Proposition 2.2 of Drappeau~\cite{drappeausommes} to prove Major Arc Estimate 1, in the appendix.

We will of course also need some information about $V(x,y;q,\delta)$.
\begin{major2}
There exist two absolute positive constants $c_{1},c_{2}$ such that the following is true. For any large $\log^{c_{1}}x \leq y \leq x$, any $q \leq y^{1/4}$, and any $|\delta| \leq \min\{y^{c_{2}}/x, e^{c_{2}\log^{1/4}x}/x\}$, we have
\begin{eqnarray}
V(x,y;q,\delta) & = & \frac{\Psi(x,y)}{x} V(x,x;q,\delta) + \nonumber \\
&& + O\left(\Psi(x,y) \frac{\log(u+1)}{\log y} \frac{2^{\omega(q)} q^{1-\alpha} \log^{2}(q+1)}{\phi(q)} \frac{\log^{3}(2+|\delta x|)}{(1+|\delta x|)^{\alpha}} \right), \nonumber
\end{eqnarray}
where $\omega(q)$ denotes the number of distinct prime factors of $q$, and $\alpha=\alpha(x,y)$ denotes the saddle-point corresponding to the $y$-smooth numbers less than $x$.
\end{major2}

When $\log^{c_{1}}x \leq y \leq e^{\sqrt{\log x}}$, say, Major Arc Estimate 2 follows from Proposition 2.3 of Drappeau~\cite{drappeausommes} (in fact with a more precise error term, and on a wider range of $q$ and $\delta$), except that the main term in Drappeau's result is expressed in terms of a certain two variable Mellin transform. It is more convenient for most applications to have a main term $(\Psi(x,y)/x)V(x,x;q,\delta)$, which will allow us simply to compare the exponential sum over smooth numbers with the complete sum, and in deducing that the error term is necessarily degraded a bit. When $y$ is larger, Drappeau's error term is less precise than in Major Arc Estimate 2 for a significant but technical reason, but we can instead deduce the result from an older estimate of La Bret\`{e}che~\cite{dlbAA}. One could presumably also use Proposition 2.4 of Drappeau~\cite{drappeausommes}, but the computations needed seem more formidable. We give details of all these deductions in the appendix.

\subsection{Other analytic tools}
In this subsection we state two general harmonic analysis results that we shall need for the proof of Theorem 2.

The first result is a very well known inequality that gives information about the (equi)distribution of a sequence in terms of bounds on its exponential sums.
\begin{harmonic1}[Erd\H{o}s--Tur\'{a}n inequality, 1948, see Corollary 1.1 of Montgomery~\cite{mont}]
Let $0 \leq u_{1},u_{2},...,u_{N} \leq 1$ be any points. Then for any $J \in \N$, and any $\alpha \leq \beta \leq \alpha + 1$, we have
$$ \left|\#\{1 \leq n \leq N : u_{n} \in [\alpha,\beta] \; \text{mod} \; 1\} - (\beta-\alpha)N \right| \leq \frac{N}{J+1} + 3\sum_{j=1}^{J} \frac{1}{j} \left|\sum_{n=1}^{N} e(ju_{n})\right| . $$
\end{harmonic1}

The second result is a tool from Bourgain's~\cite{bourgain} restriction argument, which essentially encodes the fact that major arcs are fairly ``simple'' or ``well spaced''. See Lemma 2 of Br\"{u}dern~\cite{brudern} for a closely related estimate.
\begin{harmonic2}[See pp 305-307 of Bourgain~\cite{bourgain}]
Let $x$ be large, let $Q \geq 1$ and $1/x \leq \Delta \leq 1/2$ be any parameters, and let $\epsilon > 0$ and $A > 0$ be arbitrary. Define
$$ G(\theta) := \sum_{q \leq Q} \frac{1}{q} \sum_{a=0}^{q-1} \frac{\textbf{1}_{||\theta - a/q|| \leq \Delta}}{1+x||\theta - a/q||} , $$
where as usual $||\cdot||$ denotes distance to the nearest integer, and $\textbf{1}$ denotes the indicator function.

Then if $\theta_{1},...,\theta_{R}$ are any real numbers such that $||\theta_{r} - \theta_{s}|| \geq 1/x$ when $r \neq s$,
$$ \sum_{1 \leq r,s \leq R} G(\theta_{r}-\theta_{s}) \ll_{\epsilon,A} RQ^{\epsilon}\log(1+\Delta x) + \frac{R^{2}Q\log(1+\Delta x)}{x} + \frac{R^{2}\log(1+\Delta x)}{Q^{A}} . $$
\end{harmonic2}

To orient the reader, we note that we shall apply Harmonic Analysis Result 2 towards the end of the proof of Theorem 2 (once a function like $G(\theta)$ has appeared), and in a situation where $Q$ is large compared with some other parameters in the argument, but fairly small compared with $x$. Thus the term $RQ^{\epsilon}\log(1+\Delta x)$, which basically reflects the diagonal contribution from terms $r=s$, will be the dominant one, and this nice behaviour will finally prove Theorem 2.

\section{Proof of Theorem 1}
For simplicity of writing, in this section we set $L := 2(1 + |\delta x|)$. Note that the conditions of Theorem 1 imply that $q^{2}y^{3}L^{2} \leq x$. Recall also that we write $P(m)$ for the largest prime factor of the integer $m$, and write $p(n)$ for the smallest prime factor of the integer $n$.

To prove Theorem 1, we begin with a standard move of decomposing $\sum_{\substack{n \leq x, \\ n \in \mathcal{S}(y)}} e(n\theta)$ into multiple sums, similarly as in the proof of Smooth Numbers Result 3. We see
\begin{eqnarray}
\sum_{\substack{n \leq x, \\ n \in \mathcal{S}(y)}} e(n\theta) & = & \sum_{\substack{x/qL \leq n \leq x, \\ n \in \mathcal{S}(y)}} e(n\theta) + O(\Psi(x/qL,y)) \nonumber \\
& = & \sum_{\substack{qLy < m \leq qLy^{2}, \\ m/P(m) \leq qLy, \\ m \in \mathcal{S}(y)}} \sum_{\substack{x/qLm \leq n \leq x/m, \\ p(n) \geq P(m), \\ n \in \mathcal{S}(y)}} e(mn\theta) + O(\Psi(x/qL,y)) , \nonumber
\end{eqnarray}
noting that $qLy^{2} < x/qL$, by hypothesis, so certainly every $y$-smooth number from the interval $[x/qL,x]$ has a unique decomposition in the form $mn$ (by taking $m$ to consist of the smallest prime factors of the number).

We wish to continue by applying the Cauchy--Schwarz inequality, which will allow us to complete the sum over $m$ so that we can perform estimations. However, in order to do this we need to remove the dependence between $n$ and $m$ in the condition $p(n) \geq P(m)$, and to do it efficiently we would also like the $m$ variable to run over dyadic ranges. Thus we write
\begin{eqnarray}
\sum_{\substack{n \leq x, \\ n \in \mathcal{S}(y)}} e(n\theta) & = & \sum_{0 \leq j \leq \frac{\log y}{\log 2}} \sum_{p \leq y} \sum_{\substack{2^{j}qLy < m \leq qLy \cdot \min\{2^{j+1}, p\}, \\ P(m)=p}} \sum_{\substack{x/qLm \leq n \leq x/m, \\ p(n) \geq p, \\ n \in \mathcal{S}(y)}} e(mn\theta) \nonumber \\
&& + O(\Psi(x/qL,y)) , \nonumber
\end{eqnarray}
noting that if $P(m)=p$ then the condition that $m$ be $y$-smooth is automatically satisfied, and the condition that $m/P(m) \leq qLy$ can be written as $m \leq qLyp$.

Now let $0 \leq j \leq (\log y)/\log 2$ be fixed. For simplicity of writing, we will let $\sum_{m}$ denote $\sum_{\substack{2^{j}qLy < m \leq qLy \cdot \min\{2^{j+1}, p\}, \\ P(m)=p}}$. Then using the Cauchy--Schwarz inequality, we see
\begin{eqnarray}
&& \Biggl|\sum_{p \leq y} \sum_{m} \sum_{\substack{x/qLm \leq n \leq x/m, \\ p(n) \geq p, \\ n \in \mathcal{S}(y)}} e(mn\theta) \Biggr| \leq \sqrt{\Biggl(\sum_{p \leq y} \sum_{m} 1 \Biggr) \Biggl( \sum_{p \leq y} \sum_{m} \Biggl| \sum_{\substack{x/qLm \leq n \leq x/m, \\ p(n) \geq p, \\ n \in \mathcal{S}(y)}} e(mn\theta) \Biggr|^{2} \Biggr) } \nonumber \\
& \leq & \sqrt{\Psi(2^{j+1}qLy,y)} \sqrt{\sum_{2^{j} \leq p \leq y} \sum_{\frac{2^{j}qLy}{p} < m' \leq \frac{qLy}{p} \min\{2^{j+1}, p\}} \Biggl| \sum_{\substack{x/(qLm'p) \leq n \leq x/(m'p), \\ p(n) \geq p, \\ n \in \mathcal{S}(y)}} e(m'pn\theta) \Biggr|^{2} } \nonumber \\
& \ll & \sqrt{\Psi(2^{j+1}qLy,y)} \sqrt{\sum_{2^{j} \leq p \leq y} \sum_{\substack{n_{1},n_{2} \leq \frac{x}{2^{j}qLy}, \\ p(n_{1}),p(n_{2}) \geq p, \\ n_{1}, n_{2} \in \mathcal{S}(y)}} \min\left\{\frac{2^{j+1}qLy}{p}, \frac{1}{||(n_{1}-n_{2})p\theta ||}\right\}} , \nonumber
\end{eqnarray}
where as usual $||\cdot||$ denotes distance to the nearest integer, and the final inequality follows by expanding the square and summing the geometric progression over $m'$. Notice we may restrict the sum over primes to $p \geq 2^{j}$, since otherwise the sum over $m$ (or $m'$) is empty. This small refinement is not very important, but will lead to a slightly neater and stronger bound at the end.

Now there is an unfortunate but standard complication, in that we must distinguish cases according as a prime $p$ does or does not divide $q$ (where $\theta = a/q + \delta$). Recall that $L := 2(1+|\delta x|)$. If $p \nmid q$, and if $n_{1}, n_{2} \leq x/2^{j}qLy$, then we see
$$ ||(n_{1}-n_{2})p\theta|| = || \frac{(n_{1}-n_{2})pa}{q} + (n_{1}-n_{2})p\delta || \asymp \left\{ \begin{array}{ll}
     || \frac{bpa}{q} || & \text{if} \; n_{1}-n_{2} \equiv b \; \text{mod} \; q, \; \text{for} \; b \neq 0 \\
     |(n_{1}-n_{2})p\delta| & \text{if} \; n_{1}-n_{2} \equiv 0 \; \text{mod} \; q,
\end{array} \right. $$
since by definition of $L$ we have $|(n_{1}-n_{2})p\delta| \leq p|\delta x|/2^{j}qLy \leq |\delta x|/qL \leq 1/2q$. This inequality is one of the reasons for setting the ranges of $n$ and $m$ as we did. If instead $p \mid q$, then we similarly have
$$ ||(n_{1}-n_{2})p\theta|| = || \frac{(n_{1}-n_{2})a}{q/p} + (n_{1}-n_{2})p\delta || \asymp \left\{ \begin{array}{ll}
     || \frac{ba}{q/p} || & \text{if} \; n_{1}-n_{2} \equiv b \; \text{mod} \; \frac{q}{p}, \; \text{for} \; b \neq 0 \\
     |(n_{1}-n_{2})p\delta| & \text{if} \; n_{1}-n_{2} \equiv 0 \; \text{mod} \; \frac{q}{p} .
\end{array} \right. $$

So to summarise our work so far, if we define
$$ \mathcal{T}_{j}(r) := \max_{1 \leq b \leq r-1} \sum_{\substack{n_{1},n_{2} \leq \frac{x}{2^{j}qLy}, \\ n_{1}, n_{2} \in \mathcal{S}(y), \\ n_{1}-n_{2} \equiv b \; \text{mod} \; r}} 1  , $$
(suppressing mention of the quantities $x,y,q,L$ on which $\mathcal{T}_{j}(r)$ of course also depends, and throwing away the condition that $p(n_{1}),p(n_{2}) \geq p$), then we deduce that
\begin{equation}\label{minorsummary}
\left|\sum_{\substack{n \leq x, \\ n \in \mathcal{S}(y)}} e(n\theta) \right| \ll \Psi(x/qL,y) + \sum_{0 \leq j \leq \frac{\log y}{\log 2}} \sqrt{\Psi(2^{j+1}qLy,y)} \sqrt{S_{1}} + \sum_{0 \leq j \leq \frac{\log y}{\log 2}} \sqrt{\Psi(2^{j+1}qLy,y)} \sqrt{S_{2}} ,
\end{equation}
where
$$ S_{1} = S_{1}(j) := \sum_{\substack{2^{j} \leq p \leq y, \\ p \nmid q}} \mathcal{T}_{j}(q) \sum_{b=1}^{q-1} \min\left\{\frac{2^{j+1}qLy}{p}, \frac{q}{b}\right\} + \sum_{\substack{2^{j} \leq p \leq y, \\ p \mid q}} \mathcal{T}_{j}(q/p) \sum_{b=1}^{(q/p)-1} \min\left\{\frac{2^{j+1}qLy}{p}, \frac{q}{pb}\right\} , $$
and where
\begin{eqnarray}
S_{2} = S_{2}(j) & := & \sum_{\substack{2^{j} \leq p \leq y, \\ p \nmid q}} \frac{1}{p} \sum_{\substack{n_{1},n_{2} \leq \frac{x}{2^{j}qLy}, \\ n_{1}, n_{2} \in \mathcal{S}(y), \\ n_{1}-n_{2} \equiv 0 \; \text{mod} \; q}} \min\left\{2^{j+1}qLy, \frac{1}{|(n_{1}-n_{2})\delta|}\right\} \nonumber \\
& + & \sum_{\substack{2^{j} \leq p \leq y, \\ p \mid q}} \frac{1}{p} \sum_{\substack{n_{1},n_{2} \leq \frac{x}{2^{j}qLy}, \\ n_{1}, n_{2} \in \mathcal{S}(y), \\ n_{1}-n_{2} \equiv 0 \; \text{mod} \; q/p}} \min\left\{2^{j+1}qLy, \frac{1}{|(n_{1}-n_{2})\delta|}\right\} . \nonumber
\end{eqnarray}
Here we used the fact that, if $(a,q)=1$ and $p \nmid q$, the numbers $bpa$ vary over all non-zero residue classes modulo $q$ as the numbers $b$ vary over all non-zero residue classes (though not in general in the same order, of course), similarly for the numbers $ba$ modulo $q/p$ in the case where $p \mid q$.

\vspace{12pt}
Thus far we haven't deployed much information about the $y$-smooth numbers, apart from their good factorisation properties. To complete the proof of Theorem 1 we need to bound the sums $\mathcal{T}_{j}(q), \mathcal{T}_{j}(q/p)$, and the other remaining sums over $n_{1},n_{2}$ inside $S_{2}$, and in doing this we shall use Smooth Numbers Results 2 and 3. In fact we can show:
\begin{prop1}
Let $\log^{1.1}x \leq y \leq x^{1/3}$ and $q \geq 1$ and $L = 2(1+|\delta x|)$ be as above, and suppose that $q^{2}y^{3}L^{2} \leq x$. Then for any $0 \leq j \leq (\log y)/\log 2$, and any prime $p \mid q$, we have
$$ \mathcal{T}_{j}(q) \ll \frac{\Psi(x/2^{j}qLy, y)^{2}}{q} q^{1-\alpha(x,y)} \log x \;\;\;\;\; \text{and} \;\;\;\;\; \mathcal{T}_{j}(q/p) \ll \frac{\Psi(x/2^{j}qLy, y)^{2}}{q/p} (q/p)^{1-\alpha(x,y)} \log x . $$

Under the same hypotheses, we have
$$ \sum_{\substack{n_{1},n_{2} \leq \frac{x}{2^{j}qLy}, \\ n_{1}, n_{2} \in \mathcal{S}(y), \\ n_{1}-n_{2} \equiv 0 \; \text{mod} \; q}} \min\left\{2^{j+1}qLy, \frac{1}{|(n_{1}-n_{2})\delta|}\right\} \ll 2^{j}y \Psi(\frac{x}{2^{j}qLy}, y)^{2} (qL)^{1-\alpha(x,y)} \log x \log L , $$
and similarly for any prime $p \mid q$ we have
$$ \sum_{\substack{n_{1},n_{2} \leq \frac{x}{2^{j}qLy}, \\ n_{1}, n_{2} \in \mathcal{S}(y), \\ n_{1}-n_{2} \equiv 0 \; \text{mod} \; q/p}} \min\left\{2^{j+1}qLy, \frac{1}{|(n_{1}-n_{2})\delta|} \right\} \ll p 2^{j}y \Psi(\frac{x}{2^{j}qLy}, y)^{2} (\frac{q}{p}L)^{1-\alpha(x,y)} \log x \log L . $$
\end{prop1}

Assuming the truth of Proposition 1 for a moment, we can quickly finish the proof of Theorem 1. Note that for any $p \leq y$ we have
$$ \sum_{b=1}^{q-1} \min\left\{\frac{2^{j+1}qLy}{p}, \frac{q}{b}\right\} = \sum_{b=1}^{q-1} \frac{q}{b} \ll q\log q , \;\;\; \text{and} \;\;\; \sum_{b=1}^{(q/p)-1} \min\left\{\frac{2^{j+1}qLy}{p}, \frac{q}{pb}\right\} \ll (q/p)\log q , $$
so we deduce from $(\ref{minorsummary})$ and Proposition 1 that $\left|\sum_{\substack{n \leq x, \\ n \in \mathcal{S}(y)}} e(n\theta) \right|$ is
\begin{eqnarray}
& \ll & \Psi(x/qL,y) + \sum_{0 \leq j \leq \frac{\log y}{\log 2}} \sqrt{\Psi(2^{j+1}qLy,y)} \sqrt{\frac{y}{\log y} \log q \Psi(\frac{x}{2^{j}qLy}, y)^{2} q^{1-\alpha} \log x } \nonumber \\
&& + \sum_{0 \leq j \leq \frac{\log y}{\log 2}} \sqrt{\Psi(2^{j+1}qLy,y)} \sqrt{\left(\sum_{\substack{2^{j} \leq p \leq y, \\ p \nmid q}} \frac{1}{p} + \sum_{\substack{2^{j} \leq p \leq y, \\ p \mid q}} 1 \right) 2^{j}y \Psi(\frac{x}{2^{j}qLy}, y)^{2} (qL)^{1-\alpha} \log x \log L } . \nonumber
\end{eqnarray}
Next we use the crude upper bound $\Psi(2^{j+1}qLy,y) \leq 2^{j+1}qLy$, (which, however, shouldn't lose much unless $qL$ is rather large, in which case we already gain a lot elsewhere), together with the bounds $\Psi(x/qL,y) \ll \Psi(x,y)/(qL)^{\alpha}$ and $\Psi(x/2^{j}qLy,y) \ll \Psi(x,y)/(2^{j}qLy)^{\alpha}$ coming from Smooth Numbers Result 2, and the elementary estimate
$$ \sum_{\substack{2^{j} \leq p \leq y, \\ p \nmid q}} \frac{1}{p} + \sum_{\substack{2^{j} \leq p \leq y, \\ p \mid q}} 1 \ll \sum_{p \leq y} \frac{\log p}{(j+1) p} + \sum_{\substack{p|q, \\ p \geq 2^{j}}} 1 \ll \frac{\log(qy)}{j+1} . $$
Inserting these bounds, we obtain that
\begin{eqnarray}
\left|\sum_{\substack{n \leq x, \\ n \in \mathcal{S}(y)}} e(n\theta) \right| & \ll & \frac{\Psi(x,y)}{(qL)^{\alpha}} + \sum_{0 \leq j \leq \frac{\log y}{\log 2}} \frac{\Psi(x,y)}{\sqrt{2^{j}qL}} \sqrt{\frac{1}{\log y} \log q (2^{j}qLy)^{2(1-\alpha)} q^{1-\alpha} \log x } \nonumber \\
&& + \frac{\Psi(x,y)}{\sqrt{qL}} \sum_{0 \leq j \leq \frac{\log y}{\log 2}} \sqrt{\frac{\log(qy)}{j+1} (2^{j}qLy)^{2(1-\alpha)} (qL)^{1-\alpha} \log x \log L } . \nonumber
\end{eqnarray}
In Theorem 1 we have $y \geq \log^{3}x$, and so $\alpha(x,y) \geq 2/3 + o(1)$ (by the approximation (\ref{alphaapprox})) and so the first term is more than good enough for the theorem. The second term is visibly smaller than the third, and the third term is
\begin{eqnarray}
& \ll & \frac{\Psi(x,y)}{\sqrt{qL}} (qL)^{(3/2)(1-\alpha)} \sqrt{\log(qy) \log x \log L} \sum_{0 \leq j \leq \frac{\log y}{\log 2}} \frac{2^{j(1-\alpha)} y^{1-\alpha}}{\sqrt{j+1}} \nonumber \\
& \ll & \frac{\Psi(x,y)}{\sqrt{qL}} (qL)^{(3/2)(1-\alpha)} \sqrt{\log(qy) \log x \log L} \Biggl(y^{(3/2)(1-\alpha)} \sqrt{\log y} + \sum_{\frac{\log y}{2\log 2} \leq j \leq \frac{\log y}{\log 2}} \frac{2^{j(1-\alpha)} y^{1-\alpha}}{\sqrt{\log y}} \Biggr) \nonumber \\
& \ll & \frac{\Psi(x,y)}{\sqrt{qL}} (qL)^{(3/2)(1-\alpha)} \sqrt{\log(qy) \log x \log L} \sqrt{\log y} \frac{y^{2(1-\alpha)}}{1 + |1-\alpha| \log y} . \nonumber
\end{eqnarray}
Since we have $\alpha = 1 - \frac{\log(u\log u) + O(1)}{\log y}$ and $y^{1-\alpha} = O(u\log u)$, by (\ref{alphaapprox}), the conclusion of Theorem 1 immediately follows on remembering that $u := (\log x)/\log y$.
\begin{flushright}
Q.E.D.
\end{flushright}

Now it only remains to prove Proposition 1. However, by definition we have
$$ \mathcal{T}_{j}(q) = \max_{1 \leq b \leq q-1} \sum_{\substack{n_{1} \leq x/2^{j}qLy, \\ n_{1} \in \mathcal{S}(y)}} \sum_{\substack{n_{2} \leq x/2^{j}qLy, \\ n_{2} \equiv n_{1} - b \; \text{mod} \; q, \\ n_{2} \in \mathcal{S}(y)}} 1 , $$
and by applying Smooth Numbers Result 3 on dyadic intervals $[X,2X]$, beginning with $X=qy$, we see that the inner sum is
$$ \ll \sum_{\substack{n_{2} \leq qy, \\ n_{2} \equiv n_{1} - b \; \text{mod} \; q, \\ n_{2} \in \mathcal{S}(y)}} 1 + \sum_{X = 2^{i}qy \leq x/2^{j}qLy} \frac{\Psi(X,y)}{q} q^{1-\alpha(X,y)} \log X \ll y + \frac{\Psi(x/2^{j}qLy,y)}{q} q^{1-\alpha(x,y)} \log x . $$
Here we bounded the first sum trivially, and we used Smooth Numbers Result 2 to sum the terms $\Psi(X,y)$ over all the dyadic values of $X$. We also note that
$$ \Psi(x/2^{j}qLy,y) \gg \Psi(x/2^{j}q^{2}Ly,y) q^{\alpha} \geq \Psi(Ly,y) q^{\alpha} \geq y q^{\alpha} , $$
in view of Smooth Numbers Result 2 and the fact that $x/2^{j}q^{2}Ly \geq x/q^{2}Ly^{2} \geq Ly$, by assumption. This implies that $y \ll \frac{\Psi(x/2^{j}qLy,y)}{q} q^{1-\alpha} \log x$, so we indeed have
$$ \mathcal{T}_{j}(q) \ll \sum_{\substack{n_{1} \leq x/2^{j}qLy, \\ n_{1} \in \mathcal{S}(y)}} \frac{\Psi(x/2^{j}qLy,y)}{q} q^{1-\alpha} \log x = \frac{\Psi(x/2^{j}qLy,y)^{2}}{q} q^{1-\alpha} \log x , $$
as claimed in Proposition 1. The claimed bound $\mathcal{T}_{j}(q/p) \ll \frac{\Psi(x/2^{j}qLy, y)^{2}}{q/p} (q/p)^{1-\alpha} \log x$ follows in exactly the same way.

To prove the second pair of bounds in Proposition 1, for the double sums of the form $\sum_{\substack{n_{1},n_{2} \leq \frac{x}{2^{j}qLy}, \\ n_{1}, n_{2} \in \mathcal{S}(y), \\ n_{1}-n_{2} \equiv 0 \; \text{mod} \; q}} \min\left\{2^{j+1}qLy, \frac{1}{|(n_{1}-n_{2})\delta|}\right\}$, we distinguish two cases. If $|\delta| \leq 1/x$ then we have $L = 2(1+|\delta x|) \asymp 1$, and the bounds can be proved exactly as above on upper bounding $\min\left\{2^{j+1}qLy, \frac{1}{|(n_{1}-n_{2})\delta|}\right\}$ by $2^{j+1}qLy \ll 2^{j}qy$.

The other case is where $|\delta| > 1/x$, and so we have $L \asymp |\delta x|$. For simplicity of writing, let us temporarily use $\sum^{\dagger}$ to denote a sum over pairs of integers $n_{1},n_{2} \leq \frac{x}{2^{j}qLy}$ that are $y$-smooth, and satisfy $n_{1}-n_{2} \equiv 0 \; \text{mod} \; q$. Then we have
$$ \sum_{|n_{1}-n_{2}| \leq x/2^{j}qL^{2}y}^{\dagger} \min\left\{2^{j+1}qLy, \frac{1}{|(n_{1}-n_{2})\delta|}\right\} \ll 2^{j}qLy \sum_{\substack{n_{1} \leq \frac{x}{2^{j}qLy}, \\ n_{1} \in \mathcal{S}(y)}} \sum_{\substack{|n_{2}-n_{1}| \leq x/2^{j}qL^{2}y, \\ n_{2} \in \mathcal{S}(y), \\ n_{2} \equiv n_{1} \; \text{mod} \; q}} 1 , $$
and as above we have $x/2^{j}qL^{2}y \geq qy$ by hypothesis, so Smooth Numbers Results 2 and 3 imply this is all
$$ \ll 2^{j}qLy \sum_{\substack{n_{1} \leq \frac{x}{2^{j}qLy}, \\ n_{1} \in \mathcal{S}(y)}} \frac{\Psi(x/2^{j}qLy,y)}{qL} (qL)^{1-\alpha} \log x = 2^{j}y \Psi(x/2^{j}qLy,y)^{2} (qL)^{1-\alpha} \log x . $$
Similarly, for any $0 \leq r \leq (\log L)/\log 2$ we have
\begin{eqnarray}
&& \sum_{\frac{2^{r}x}{2^{j}qL^{2}y} < |n_{1}-n_{2}| \leq \frac{2^{r+1}x}{2^{j}qL^{2}y}}^{\dagger} \min\left\{2^{j+1}qLy, \frac{1}{|(n_{1}-n_{2})\delta|}\right\} \ll \frac{2^{j}qLy}{2^{r}} \sum_{\substack{n_{1} \leq \frac{x}{2^{j}qLy}, \\ n_{1} \in \mathcal{S}(y)}} \sum_{\substack{|n_{2}-n_{1}| \leq \frac{2^{r+1}x}{2^{j}qL^{2}y}, \\ n_{2} \in \mathcal{S}(y), \\ n_{2} \equiv n_{1} \; \text{mod} \; q}} 1 \nonumber \\
& \ll & \frac{2^{j}qLy}{2^{r}} \sum_{\substack{n_{1} \leq \frac{x}{2^{j}qLy}, \\ n_{1} \in \mathcal{S}(y)}} \Psi(x/2^{j}qLy,y) \left(\frac{2^{r}}{qL}\right)^{\alpha} \log x \nonumber \\
& \ll & 2^{j}y \Psi(x/2^{j}qLy,y)^{2} (qL)^{1-\alpha} \log x , \nonumber
\end{eqnarray}
and so the claimed bound for $\sum^{\dagger} \min\left\{2^{j+1}qLy, \frac{1}{|(n_{1}-n_{2})\delta|}\right\}$ follows by summing over $r$.

The claimed bound for $\sum_{\substack{n_{1},n_{2} \leq \frac{x}{2^{j}qLy}, \\ n_{1}, n_{2} \in \mathcal{S}(y), \\ n_{1}-n_{2} \equiv 0 \; \text{mod} \; q/p}} \min\left\{2^{j+1}qLy, \frac{1}{|(n_{1}-n_{2})\delta|} \right\}$, losing a factor $p^{\alpha}$ because of the weaker congruence condition, follows in exactly the same way.
\begin{flushright}
Q.E.D.
\end{flushright}

\section{Proof of Theorem 2}
Throughout this section we shall assume that $y \leq x^{1/100}$, since if $x^{1/100} < y \leq x$ then $\Psi(x,y) \asymp x$, in which case Theorem 2 is a trivial consequence of Parseval's identity.

To prove Theorem 2, we shall actually prove the following large values (or ``distributional'') estimate for the number of $1/x$-neighbourhoods on which $\left| \sum_{\substack{n \leq x, \\ n \in \mathcal{S}(y)}} a_{n} e(n\theta) \right|$ is large. Recall that $||\cdot||$ denotes distance to the nearest integer.

\begin{prop2}
There exist large constants $C_{0}, c_{1} > 0$ such that the following is true.

Let $\log^{21}x \leq y \leq x^{1/100}$ be large, and let $(a_{n})_{n \leq x}$ be any complex numbers with absolute values at most 1. Let $x^{-1/1000} < \delta \leq 1$, and let $0 \leq \theta_{1}, \theta_{2}, ..., \theta_{R} \leq 1$ be any real numbers such that $||\theta_{r} - \theta_{s}|| \geq 1/x$ when $r \neq s$, and such that
$$ \Biggl| \sum_{\substack{n \leq x, \\ n \in \mathcal{S}(y)}} a_{n} e(n\theta_{r}) \Biggr| \geq \delta \Psi(x,y) \;\;\;\;\; \forall 1 \leq r \leq R. $$

Then for any small $\epsilon > 0$ we have
\begin{equation}\label{smalldelta}
R \ll_{\epsilon} \delta^{-2 - 20(1-\alpha(x,y)) - \epsilon} \log^{3+32(1-\alpha(x,y))+\epsilon}x ,
\end{equation}
where $\alpha(x,y)$ denotes the saddle-point corresponding to the $y$-smooth numbers below $x$.

In addition, if $\log^{c_{1}}x \leq y \leq x^{1/100}$ and $\delta^{-1} \leq \min\{\log^{1/\epsilon}x, y^{1/C_{0}}\}$ then
\begin{equation}\label{largedelta}
R \ll_{\epsilon} \delta^{-2 - 8(1-\alpha(x,y)) - \epsilon} .
\end{equation}
\end{prop2}

First let us use Proposition 2 to prove Theorem 2. Recall that we have $1-\alpha = \frac{\log(u\log(u+1)) + O(1)}{\log y} \leq 1/20$, say, in view of the approximation (\ref{alphaapprox}) and our assumption that $y \geq \log^{21}x$. If $\delta^{-1} > \log^{1/\epsilon}x$ then
$$ \log^{3+32(1-\alpha)+\epsilon}x \leq \log^{5}x \leq \delta^{-5\epsilon}, $$
or if $y^{1/C_{0}} < \delta^{-1} \leq \log^{1/\epsilon}x$ then $u := (\log x)/\log y \gg_{\epsilon} (\log x)/\log\log x$, and therefore
$$ \log^{3+32(1-\alpha)+\epsilon}x \leq \log^{5}x \ll_{\epsilon} u^{6} \ll y^{6(1-\alpha)} < \delta^{-6C_{0}(1-\alpha)} . $$
So after relabelling the arbitrary small $\epsilon$, we can combine (\ref{smalldelta}) and (\ref{largedelta}) into the single bound
$$ R \ll_{\epsilon} \delta^{-2-(20+6C_{0})(1-\alpha)-\epsilon} . $$
Applying this bound, together with Parseval's identity, we deduce that
\begin{eqnarray}
&& \int_{0}^{1} \left| \sum_{\substack{n \leq x, \\ n \in \mathcal{S}(y)}} a_{n} e(n\theta) \right|^{p} d\theta \nonumber \\
& \leq & \sum_{1 \leq j \leq \frac{\log x}{1000\log 2}} \left(\frac{\Psi(x,y)}{2^{j-1}}\right)^{p} \text{meas}\Biggl\{0 \leq \theta \leq 1 : \frac{\Psi(x,y)}{2^{j}} < \Biggl| \sum_{\substack{n \leq x, \\ n \in \mathcal{S}(y)}} a_{n} e(n\theta) \Biggr| \leq \frac{\Psi(x,y)}{2^{j-1}} \Biggr\} \nonumber \\
&& + \left(\frac{\Psi(x,y)}{x^{1/1000}/2}\right)^{p-2} \int_{0}^{1} \Biggl| \sum_{\substack{n \leq x, \\ n \in \mathcal{S}(y)}} a_{n} e(n\theta) \Biggr|^{2} d\theta \nonumber \\
& \ll_{p, \epsilon} & \sum_{1 \leq j \leq \frac{\log x}{1000\log 2}} \left(\frac{\Psi(x,y)}{2^{j}}\right)^{p} \cdot \frac{1}{x} 2^{j(2 + (20+6C_{0})(1-\alpha) + \epsilon)} + \frac{\Psi(x,y)^{p-1}}{x^{(p-2)/1000}} . \nonumber
\end{eqnarray}
We assume in Theorem 2 that $y \geq \log^{C\max\{1,1/(p-2)\}}x$, which implies by (\ref{alphaapprox}) that
$$ 1 - \alpha(x,y) \leq \frac{\log\log x}{\log y} + o(1) \leq \frac{\min\{1,p-2\}}{C} + o(1) , $$
 and which also implies that $\Psi(x,y)/x \geq \Psi(x,\log^{C\max\{1,1/(p-2)\}}x)/x \geq x^{-\min\{1,p-2\}/C + o(1)}$. So if we choose $\epsilon=(p-2)/3$, say, then provided $C$ is large enough we will have
$$ \sum_{1 \leq j \leq \frac{\log x}{1000\log 2}} \left(\frac{\Psi(x,y)}{2^{j}}\right)^{p} \cdot \frac{1}{x} 2^{j(2 + (20+6C_{0})(1-\alpha) + \epsilon)} \leq \frac{\Psi(x,y)^{p}}{x} \sum_{1 \leq j \leq \frac{\log x}{1000\log 2}} \frac{1}{2^{j(p-2)/3}} \ll_{p} \frac{\Psi(x,y)^{p}}{x}, $$
and also $\Psi(x,y)^{p-1}/x^{(p-2)/1000} \ll \Psi(x,y)^{p}/x$. Theorem 2 follows immediately.
\begin{flushright}
Q.E.D.
\end{flushright}

It remains to prove Proposition 2, and we shall assume throughout, without loss of generality, that $\delta$ is at most a small absolute constant. We shall also concentrate on proving the bound (\ref{smalldelta}). The bound (\ref{largedelta}), which is only relevant to save logarithmic factors when $\delta$ is fairly large, can be proved by inserting Major Arc Estimates 1 and 2 into Bourgain's~\cite{bourgain} original restriction argument, as we shall briefly explain at the end.

To try to make things more digestible, we split the proof of (\ref{smalldelta}) into two halves. Firstly we shall prove the following intermediate statement, which already encodes the Bourgain/Hal\'{a}sz--Montgomery duality step, our combination of this with double sums, and the use of Theorem 1 and the Erd\H{o}s--Tur\'{a}n inequality to reduce to looking at differences $\theta_{r}-\theta_{s}$ on ``fairly major'' arcs.
\begin{prop3}
Let the situation be as in Proposition 2, with $\delta$ at most a small absolute constant, and define $K := \lfloor (\delta^{-1}\log x)^{10} \rfloor$ (where $\lfloor \cdot \rfloor$ denotes integer part).

Then
$$ \delta^{2} R^{2} \Psi(x,y) \ll \sum_{q \leq K\log^{12}x} \sum_{(a,q)=1} \sum_{\substack{1 \leq r,s \leq R, \\ ||\theta_{r}-\theta_{s}-\frac{a}{q}|| \leq \frac{K^{2} y \log^{12}x}{x}}} \sum_{\substack{x/yK < m \leq x/K, \\ m \in \mathcal{S}(y)}} \min\{\frac{x}{m},\frac{1}{||m(\theta_{r}-\theta_{s})||}\} . $$
\end{prop3}

To prove Proposition 3, let $\theta_{1},...,\theta_{R}$ be as in the statement, and let $c_{1},...,c_{R}$ be complex numbers of absolute value 1 such that $c_{r} \sum_{\substack{n \leq x, \\ n \in \mathcal{S}(y)}} a_{n} e(n\theta_{r}) = \left|\sum_{\substack{n \leq x, \\ n \in \mathcal{S}(y)}} a_{n} e(n\theta_{r}) \right|$. Then the Cauchy--Schwarz inequality yields that
\begin{eqnarray}
\delta^{2} R^{2} \Psi(x,y)^{2} \leq \left(\sum_{r=1}^{R} \left|\sum_{\substack{n \leq x, \\ n \in \mathcal{S}(y)}} a_{n} e(n\theta_{r}) \right| \right)^{2} & = & \left(\sum_{r=1}^{R} c_{r} \sum_{\substack{n \leq x, \\ n \in \mathcal{S}(y)}} a_{n} e(n\theta_{r}) \right)^{2} \nonumber \\
& \leq & \Psi(x,y) \sum_{\substack{n \leq x, \\ n \in \mathcal{S}(y)}} \left|\sum_{r=1}^{R} c_{r}e(n\theta_{r}) \right|^{2} . \nonumber
\end{eqnarray}
This manipulation follows the first part of the argument in $\S 4$ of Bourgain~\cite{bourgain}, and multiplicative number theorists may be most familiar with it from the work of Hal\'{a}sz and others on large values of Dirichlet polynomials.

At this point, Bourgain~\cite{bourgain} proceeds by expanding the square, switching the order of summation, and inserting his good major arc estimates. But we do not have sufficient major arc information, so we shall proceed differently. Note that $1 < x/yK < \delta^{3}x$, because of our assumptions that $y \leq x^{1/100}$ and $\delta > x^{-1/1000}$. Thus, if we again write $P(m)$ for the largest prime factor of the integer $m$, and write $p(n)$ for the smallest prime factor of the integer $n$, we can rewrite
\begin{eqnarray}
\sum_{\substack{n \leq x, \\ n \in \mathcal{S}(y)}} \left|\sum_{r=1}^{R} c_{r}e(n\theta_{r}) \right|^{2} & = & \sum_{\substack{x/yK < m \leq x/K, \\ m/P(m) \leq x/yK, \\ m \in \mathcal{S}(y)}} \sum_{\substack{\delta^{3}x/m \leq n \leq x/m, \\ p(n) \geq P(m), \\ n \in \mathcal{S}(y)}} \left|\sum_{r=1}^{R} c_{r}e(mn\theta_{r}) \right|^{2} + \sum_{\substack{n < \delta^{3}x, \\ n \in \mathcal{S}(y)}} \left|\sum_{r=1}^{R} c_{r}e(n\theta_{r}) \right|^{2} \nonumber \\
& \leq & \sum_{\substack{x/yK < m \leq x/K, \\ m/P(m) \leq x/yK, \\ m \in \mathcal{S}(y)}} \sum_{\delta^{3}x/m \leq n \leq x/m} \left|\sum_{r=1}^{R} c_{r}e(mn\theta_{r}) \right|^{2} + R^{2} \Psi(\delta^{3}x,y) \nonumber \\
& \ll & \sum_{\substack{x/yK < m \leq x/K, \\ m/P(m) \leq x/yK, \\ m \in \mathcal{S}(y)}} \sum_{1 \leq r,s \leq R} \min\{\frac{x}{m},\frac{1}{||m(\theta_{r}-\theta_{s})||}\} + R^{2}\Psi(x,y) \delta^{3\alpha} . \nonumber
\end{eqnarray}
Here the final inequality follows by expanding the square and performing the summation over $n$, and by using Smooth Numbers Result 2. In particular, since $y \geq \log^{21}x$ we have $\alpha(x,y) = 1 - \frac{\log(u\log(u+1)) + O(1)}{\log y} \geq 19/20$ (say), by the approximation (\ref{alphaapprox}). So provided $\delta$ is sufficiently small the contribution from the final term will certainly be $\leq \delta^{2} R^{2} \Psi(x,y)/2$, and we will have
$$ \delta^{2} R^{2} \Psi(x,y) \ll \sum_{\substack{x/yK < m \leq x/K, \\ m \in \mathcal{S}(y)}} \sum_{1 \leq r,s \leq R} \min\{\frac{x}{m},\frac{1}{||m(\theta_{r}-\theta_{s})||}\} . $$

To finish the proof of Proposition 3, we need to show that the only pairs $r,s$ that can make a substantial contribution to the right hand side are those where $\theta_{r}-\theta_{s}$ is ``fairly close'' to a rational with ``fairly small'' denominator (in terms of $K := \lfloor (\delta^{-1}\log x)^{10} \rfloor$). One could presumably show this directly by combinatorial arguments, but this would entail a further careful decomposition of the $m$ variable. Instead we shall use the Erd\H{o}s--Tur\'{a}n inequality (Harmonic Analysis Result 1), together with Theorem 1. Indeed, for any $r,s \leq R$ we have
\begin{eqnarray}
\sum_{\substack{\frac{x}{yK} < m \leq \frac{x}{K}, \\ m \in \mathcal{S}(y)}} \min\{\frac{x}{m},\frac{1}{||m(\theta_{r}-\theta_{s})||}\} & \leq & \Psi(\frac{x}{K},y) \sqrt{K} + \sum_{\substack{x/yK < m \leq x/K, \\ m \in \mathcal{S}(y), \\ ||m(\theta_{r}-\theta_{s})|| \leq 1/\sqrt{K}}} \min\{\frac{x}{m},\frac{1}{||m(\theta_{r}-\theta_{s})||}\} \nonumber \\
& \ll & \Psi(x,y) K^{1/2-\alpha} + x \sum_{\substack{x/yK < m \leq x/K, \\ m \in \mathcal{S}(y), \\ ||m(\theta_{r}-\theta_{s})|| \leq 1/\sqrt{K}}} \frac{1}{m} , \nonumber
\end{eqnarray}
where the second inequality used Smooth Numbers Result 2. For simplicity of writing, for any $x/yK < M \leq x/K$ let us temporarily write $\sum_{m \sim M}$ to mean $\sum_{\substack{M < m \leq 2M, \\ m \in \mathcal{S}(y), \\ ||m(\theta_{r}-\theta_{s})|| \leq 1/\sqrt{K}}}$. Then the Erd\H{o}s--Tur\'{a}n inequality implies that, for any $J \in \N$,
$$ \sum_{m \sim M} 1 \ll \frac{\Psi(2M,y)}{\sqrt{K}} + \frac{\Psi(2M,y)}{J+1} + \sum_{j=1}^{J} \frac{1}{j} \left|\sum_{\substack{M < m \leq 2M, \\ m \in \mathcal{S}(y)}} e(jm(\theta_{r}-\theta_{s})) \right| . $$
In particular, if we choose $J := \lfloor \sqrt{K} \rfloor$, and if $\theta_{r}-\theta_{s}=a/q+\eta$ for some $q \leq x^{0.6}$, some $(a,q)=1$ and some $|\eta| \leq 1/qx^{0.6}$ (say), then
$$ \sum_{m \sim M} 1 \ll \frac{\Psi(2M,y)}{\sqrt{K}} + \sum_{h|q} \frac{1}{h} \sum_{\substack{j' \leq \sqrt{K}/h, \\ (j',q/h)=1}} \frac{1}{j'} \left|\sum_{\substack{M < m \leq 2M, \\ m \in \mathcal{S}(y)}} e(mj'h(\theta_{r}-\theta_{s})) \right| , $$
and we will have $j'h(\theta_{r}-\theta_{s}) = aj'/(q/h) + j'h\eta$, where $aj'$ is coprime to $q/h$ and where $|j'h\eta| \leq \sqrt{K}/qx^{0.6} \leq 1/qx^{0.55}$.

Now in the case where $(q/h) \leq x^{0.45}$, say, the condition $(q/h)^{2}y^{3}(1+|j'h\eta M|)^{2} \leq M/4$ in Theorem 1 will be satisfied (remembering that $y \leq x^{1/100}$ throughout this section), and the theorem implies
$$ \frac{1}{h} \Biggl|\sum_{\substack{M < m \leq 2M, \\ m \in \mathcal{S}(y)}} e(mj'h(\theta_{r}-\theta_{s})) \Biggr| \ll \frac{1}{h} \frac{\Psi(2M,y) \log^{7/2}x}{((q/h)(1+|j'h\eta M|))^{1/2 - 3/2(1-\alpha)}} \ll \frac{\Psi(2M,y) \log^{7/2}x}{(q(1+|\eta M|))^{1/2 - 3/2(1-\alpha)}} . $$
If $x^{0.45} < (q/h) \leq x^{0.6}$ then we can instead apply Fouvry and Tenenbaum's~\cite{fouvryten0} minor arc estimate (as stated in our survey in the introduction), obtaining a bound $\ll M/x^{0.15}$ (say) which is negligible compared with the term $\Psi(2M,y)/\sqrt{K}$. We deduce overall that
\begin{eqnarray}
\sum_{m \sim M} 1 & \ll & \frac{\Psi(2M,y)}{\sqrt{K}} + \frac{\Psi(2M,y) d(q) \log^{9/2}x}{(q(1+|\eta M|))^{1/2 - 3/2(1-\alpha)}} \ll \frac{\Psi(2M,y)}{\sqrt{K}} + \frac{\Psi(2M,y) \log^{9/2}x}{(q(1+|\eta x/yK|))^{2/5}} , \nonumber
\end{eqnarray}
where $d(q)$ denotes the divisor function. Here we again used the fact that $1-\alpha = \frac{\log(u\log(u+1)) + O(1)}{\log y} \leq 1/20$, and the fact that $M > x/yK$.

In particular, if either $q \geq K \log^{12}x$ or $|\eta| \geq K^{2} y (\log^{12}x)/x$ then the above is $\ll \Psi(2M,y)/K^{2/5}$. By summing over dyadic values of $M$ it finally follows, using Smooth Numbers Result 2, that
$$ x \sum_{\substack{x/yK < m \leq x/K, \\ m \in \mathcal{S}(y), \\ ||m(\theta_{r}-\theta_{s})|| \leq 1/\sqrt{K}}} \frac{1}{m} \ll \frac{\Psi(x,y)}{K^{2/5}} \log y (yK)^{1-\alpha} = \Psi(x,y) K^{3/5-\alpha} y^{1-\alpha} \log y . $$
Then we have $K^{3/5-\alpha} \leq K^{3/5-19/20} = K^{-7/20} \ll \delta^{7/2}/\log^{7/2}x$, by definition of $K := \lfloor (\delta^{-1}\log x)^{10} \rfloor$. We also have $y^{1-\alpha} \log y \ll u \log(u+1) \log y \ll \log^{2}x$, so we see that the contribution from $\theta_{r}-\theta_{s}$ is indeed negligible compared with $\delta^{2} \Psi(x,y)$ if either $q \geq K \log^{12}x$ or $|\eta| \geq K^{2} y (\log^{12}x)/x$, and Proposition 3 follows immediately.
\begin{flushright}
Q.E.D.
\end{flushright}

To complete the proof of (\ref{smalldelta}), we need to bound $\sum_{\substack{x/yK < m \leq x/K, \\ m \in \mathcal{S}(y)}} \min\{\frac{x}{m},\frac{1}{||m(\theta_{r}-\theta_{s})||}\}$ when $||\theta_{r}-\theta_{s}-\frac{a}{q}|| \leq \frac{K^{2} y \log^{12}x}{x}$ for some $q \leq K\log^{12}x$, as in the conclusion of Proposition 3. In fact we shall prove the following result.
\begin{prop4}
Let the situation be as in Proposition 3. Suppose that $q \leq K\log^{12}x$, that $(a,q)=1$, and that $|\eta| \leq K^{2}y(\log^{12}x)/x$. Also define $\widetilde{M} := \min\{1/2q|\eta|,x/K\}$.

Then if $\theta = a/q + \eta$, we have the bounds
$$ \sum_{\substack{\widetilde{M} < m \leq x/K, \\ m \in \mathcal{S}(y)}} \min\{\frac{x}{m},\frac{1}{||m\theta||}\} \ll \frac{\delta^{2}\Psi(x,y)}{\log x} $$
and
$$ \sum_{\substack{x/yK < m \leq \widetilde{M}, \\ m \in \mathcal{S}(y)}} \min\{\frac{x}{m},\frac{1}{||m\theta||}\} \ll \frac{\delta^{2}\Psi(x,y)}{\log x} + \frac{\Psi(x,y)}{q(1+|\eta x|)} \log^{2}x (K^{2}\log^{12}x)^{1-\alpha(x,y)} . $$
\end{prop4}

The proof of Proposition 4 depends on Smooth Numbers Results 2 and 3, similarly as the proof of Theorem 1.

We begin with the first bound in Proposition 4, which is the harder one. We may assume that $|\eta| \geq K/2qx$ and therefore $\widetilde{M} = 1/2q|\eta|$, since otherwise the range of summation in the bound is empty and the statement is trivial. We may also assume without loss of generality (on negating $\theta$) that $\eta > 0$. Observe then that if $m$ lies in an interval of the form $(j/qK\eta,(j+1)/qK\eta]$ then we have $ |m\theta - ma/q - j/qK| \leq 1/qK$. Further, if $j \equiv k$ modulo $K$, where we choose the representative $k$ such that $-K/2 < k \leq K/2$, then for a certain integer $f(j)$ we will have
$$ |m\theta - ma/q - f(j)/q - k/qK| \leq 1/qK . $$
Therefore as $m$ ranges over the interval and over all residue classes modulo $q$, except for the unique class such that $ma \equiv -f(j)$ modulo $q$, the numbers $||m\theta|| \asymp ||ma/q + f(j)/q||$ will be $\asymp b/q$ for those $1 \leq b \leq q/2$. Meanwhile, if $|k| \geq 2$ then, when $m$ is in the interval and $ma \equiv -f(j)$ modulo $q$, we will have $ ||m\theta|| \asymp k/qK $. It is only if $ma \equiv -f(j)$ modulo $q$ and $k \in \{-1,0,1\}$ that we cannot usefully lower bound $||m\theta||$.

In view of the above discussion, we conclude that
\begin{eqnarray}\label{combinatoricshc}
&& \sum_{\substack{\widetilde{M} < m \leq x/K, \\ m \in \mathcal{S}(y)}} \min\{\frac{x}{m},\frac{1}{||m\theta||}\} \leq \sum_{-K/2 < k \leq K/2} \sum_{\substack{qK\widetilde{M}\eta - 1 \leq j \leq qx\eta, \\ j \equiv k \; \text{mod} \; K}} \sum_{b=1}^{q} \sum_{\substack{\frac{j}{qK\eta} < m \leq \frac{j+1}{qK\eta}, \\ m \equiv b \; \text{mod} \; q, \\ m \in \mathcal{S}(y)}} \min\{\frac{x}{m},\frac{1}{||m\theta||}\} \nonumber \\
& \ll & \sum_{|k| \leq \frac{K}{2}} \sum_{\substack{j \leq qx\eta, \\ j \equiv k \; \text{mod} \; K}} \left(\sum_{b=1}^{q} \frac{q}{b} \right) \max_{1 \leq b \leq q} \sum_{\substack{\frac{j}{qK\eta} < m \leq \frac{j+1}{qK\eta}, \\ m \equiv b \; \text{mod} \; q, \\ m \in \mathcal{S}(y)}} 1 + \sum_{\substack{|k| \leq \frac{K}{2}, \\ k \notin \{-1,0,1\}}} \frac{qK}{|k|} \sum_{\substack{j \leq qx\eta, \\ j \equiv k \; \text{mod} \; K}} \sum_{\substack{\frac{j}{qK\eta} < m \leq \frac{j+1}{qK\eta}, \\ m \equiv -a^{-1}f(j) \; \text{mod} \; q, \\ m \in \mathcal{S}(y)}} 1 \nonumber \\
&& + \sum_{-1 \leq k \leq 1} \sum_{\substack{qK\widetilde{M}\eta - 1 \leq j \leq qx\eta, \\ j \equiv k \; \text{mod} \; K}} \sum_{\substack{\frac{j}{qK\eta} < m \leq \frac{j+1}{qK\eta}, \\ m \equiv -a^{-1}f(j) \; \text{mod} \; q, \\ m \in \mathcal{S}(y)}} \frac{x}{m} .
\end{eqnarray}
Each of the first two terms in (\ref{combinatoricshc}) is
\begin{eqnarray}
& \ll & qK\log(q+1)  \sum_{|k| \leq \frac{K}{2}} \frac{1}{|k|+1} \sum_{\substack{j \leq qx\eta, \\ j \equiv k \; \text{mod} \; K}} \max_{1 \leq b \leq q} \sum_{\substack{\frac{j}{qK\eta} < m \leq \frac{j+1}{qK\eta}, \\ m \equiv b \; \text{mod} \; q, \\ m \; \text{is} \; y-\text{smooth}}} 1 \nonumber \\
& \ll & qK\log(q+1)  \sum_{|k| \leq \frac{K}{2}} \frac{1}{|k|+1} \sum_{\substack{j \leq qx\eta, \\ j \equiv k \; \text{mod} \; K}} \frac{\Psi(j/qK\eta,y)}{(jq)^{\alpha}} \log x , \nonumber
\end{eqnarray}
by Smooth Numbers Result 3. Further, using Smooth Numbers Result 2 to bound $\Psi(j/qK\eta,y)$ we find the above is
$$ \ll qK\log(q+1)  \sum_{|k| \leq \frac{K}{2}} \frac{1}{|k|+1} \sum_{\substack{j \leq qx\eta, \\ j \equiv k \; \text{mod} \; K}} \frac{\Psi(x,y)}{q^{2}K\eta x} (q^{2}K\eta x)^{1-\alpha} \log x , $$
Since we assume that $\eta \geq K/2qx$ we have $qx\eta \geq K/2$, and so the above is
$$ \ll qK\log(q+1)  \sum_{|k| \leq \frac{K}{2}} \frac{1}{|k|+1} \frac{\Psi(x,y)}{q K^{2}} (q^{2}K\eta x)^{1-\alpha} \log x \ll \frac{\Psi(x,y)}{K} \log^{3}x (q^{2}K\eta x)^{1-\alpha} . $$ 
Similarly, using the bound $x/m < xqK\eta/j$ and the Smooth Numbers Results, the third term in (\ref{combinatoricshc}) is seen to be
\begin{eqnarray}
& \ll & xqK\eta \sum_{-1 \leq k \leq 1} \sum_{\substack{qK\widetilde{M}\eta - 1 \leq j \leq qx\eta, \\ j \equiv k \; \text{mod} \; K}} \frac{1}{j} \max_{1 \leq b \leq q} \sum_{\substack{\frac{j}{qK\eta} < m \leq \frac{j+1}{qK\eta}, \\ m \equiv b \; \text{mod} \; q, \\ m \in \mathcal{S}(y)}} 1 \nonumber \\
& \ll & \frac{\Psi(x,y)}{q} (q^{2}K\eta x)^{1-\alpha} \log x \sum_{-1 \leq k \leq 1} \sum_{\substack{qK\widetilde{M}\eta - 1 \leq j \leq qx\eta, \\ j \equiv k \; \text{mod} \; K}} \frac{1}{j} \ll \frac{\Psi(x,y)}{qK} (q^{2}K\eta x)^{1-\alpha} \log^{2}x . \nonumber
\end{eqnarray}
Here the final inequality crucially uses the fact that $qK\widetilde{M}\eta = K/2$, and therefore the smallest terms in the sum over $j$ are $\asymp K$.

Since $q^{2}K\eta x \leq q^{2}K^{3}y\log^{12}x \leq K^{5}y\log^{36}x$, both of the above bounds are
$$ \ll \frac{\Psi(x,y)}{K} \log^{3}x (K^{5}\log^{36}x)^{1-\alpha} y^{1-\alpha} . $$
Moreover, we have $y^{1-\alpha} \ll u\log(u+1) \ll \log x$ and $1-\alpha \leq 1/20$ by the approximation (\ref{alphaapprox}) (and the fact that $y \geq \log^{21}x$), so the bounds are
$$ \ll \frac{\Psi(x,y)}{K^{3/4}} \log^{29/5}x . $$
Since $K \gg (\delta^{-1}\log x)^{10}$ this is all $\ll \Psi(x,y) \delta^{30/4}/\log^{17/10}x \ll \Psi(x,y)\delta^{2}/\log x$, which is precisely the first bound claimed in Proposition 4.

It remains to prove the second bound in Proposition 4, for $\sum_{\substack{x/yK < m \leq \widetilde{M}, \\ m \in \mathcal{S}(y)}} \min\{\frac{x}{m},\frac{1}{||m\theta||}\}$. This time we don't need to split the range of $m$ into intervals, since we have $|m\eta| \leq \widetilde{M}|\eta| \leq 1/2q$ and therefore we simply have
$$ ||m\theta|| = || \frac{ma}{q} + m\eta || \asymp \left\{ \begin{array}{ll}
     || \frac{ba}{q} || & \text{if} \; m \equiv b \; \text{mod} \; q, \; \text{for} \; b \neq 0 \\
     |m\eta| & \text{if} \; m \equiv 0 \; \text{mod} \; q.
\end{array} \right. $$
So we find, similarly as above (and remembering that $\widetilde{M} \leq x/K$), that
\begin{eqnarray}
\sum_{\substack{x/yK < m \leq \widetilde{M}, \\ m \in \mathcal{S}(y)}} \min\{\frac{x}{m},\frac{1}{||m\theta||}\} & \ll & \left(\sum_{b=1}^{q-1} \frac{q}{b} \right) \max_{1 \leq b \leq q-1} \sum_{\substack{x/yK < m \leq \widetilde{M}, \\ m \equiv b \; \text{mod} \; q, \\ m \in \mathcal{S}(y)}} 1 + \sum_{\substack{x/yK < m \leq \widetilde{M}, \\ m \equiv 0 \; \text{mod} \; q, \\ m \in \mathcal{S}(y)}} \min\{\frac{x}{m},\frac{1}{m|\eta|}\} \nonumber \\
& \ll & q\log q \frac{\Psi(\frac{x}{K},y)}{q} q^{1-\alpha} \log x + \frac{1}{q} \min\{x,\frac{1}{|\eta|}\} \sum_{\substack{x/qyK < m' \leq x/qK, \\ m' \in \mathcal{S}(y)}} \frac{1}{m'} . \nonumber
\end{eqnarray}
Then Smooth Numbers Result 2 and the fact that $\alpha \geq 19/20$ imply that $\Psi(x/K,y) \ll \Psi(x,y)/K^{\alpha} \ll \Psi(x,y)/K^{19/20}$, and so the first term is $\ll \Psi(x,y) q^{1/20} (\log^{2}x)/K^{19/20}$. Since we also have $q \leq K\log^{12}x$, with $K=\lfloor (\delta^{-1}\log x)^{10} \rfloor$, the first term is $\ll \Psi(x,y) (\log^{13/5}x)/K^{9/10} \ll \Psi(x,y)\delta^{2}/\log x$, which is acceptable for Proposition 4. Meanwhile, by dividing into dyadic intervals and using Smooth Numbers Result 2 we find
\begin{eqnarray}
\sum_{\substack{x/qyK < m' \leq x/qK, \\ m' \in \mathcal{S}(y)}} \frac{1}{m'} \ll \log y \max_{x/qyK < M \leq x/qK} \frac{\Psi(M,y)}{M} & \ll & \frac{\Psi(x,y)}{x} \log y (qyK)^{1-\alpha} \nonumber \\
& \ll & \frac{\Psi(x,y)}{x} \log^{2}x (qK)^{1-\alpha} , \nonumber
\end{eqnarray}
using again the fact that $y^{1-\alpha} \ll u\log(u+1) = (\log x \log(u+1))/\log y$. Since $q \leq K\log^{12}x$, we finally conclude that
$$ \sum_{\substack{x/yK < m \leq \widetilde{M}, \\ m \in \mathcal{S}(y)}} \min\{\frac{x}{m},\frac{1}{||m\theta||}\} \ll \frac{\delta^{2}\Psi(x,y)}{\log x} + \frac{\Psi(x,y)}{q(1+|\eta x|)} \log^{2}x (K^{2}\log^{12}x)^{1-\alpha} , $$
as claimed in Proposition 4.
\begin{flushright}
Q.E.D.
\end{flushright}

Finally, by combining Propositions 3 and 4 we obtain that
$$ \delta^{2}R^{2}\Psi(x,y) \ll \sum_{q \leq K\log^{12}x} \sum_{(a,q)=1} \sum_{1 \leq r,s \leq R} \frac{\Psi(x,y)}{q(1+x||\theta_{r}-\theta_{s}-a/q||)} \log^{2}x (K^{2}\log^{12}x)^{1-\alpha} . $$
Here we noted that there is at most one pair $a,q$ such that $||\theta_{r}-\theta_{s}-\frac{a}{q}|| \leq \frac{K^{2} y \log^{12}x}{x}$, which means that the contribution from all terms $\frac{\delta^{2}\Psi(x,y)}{\log x}$ in Proposition 4 to Proposition 3 is $\ll \frac{\delta^{2}R^{2}\Psi(x,y)}{\log x}$, which is negligible.

If we now apply Harmonic Analysis Result 2, with the choices $Q=K\log^{12}x$ and $\Delta = 1/2$ and $A=1$, we deduce that
$$ \delta^{2}R^{2} \ll_{\epsilon} \log^{2}x (K^{2}\log^{12}x)^{1-\alpha} \left(RQ^{\epsilon}\log x + \frac{R^{2}Q\log x}{x} + \frac{R^{2}\log x}{Q} \right) . $$
Since $1-\alpha \leq 1/20$ and $K=\lfloor (\delta^{-1}\log x)^{10} \rfloor$, and $\delta \geq x^{-1/1000}$ and $Q \leq x^{0.1}$ (say), the contribution from the last two terms is negligible and we obtain
$$ R \ll_{\epsilon} \log^{3}x (K^{2}\log^{12}x)^{1-\alpha} \delta^{-2} Q^{\epsilon} \ll_{\epsilon} \delta^{-2 - 20(1-\alpha) - 10\epsilon} \log^{3+32(1-\alpha)+22\epsilon}x , $$
which is precisely the claimed bound (\ref{smalldelta}) after relabelling the arbitrary small $\epsilon$.
\begin{flushright}
Q.E.D.
\end{flushright}

As promised, we end this section by briefly sketching how to obtain the other bound (\ref{largedelta}). Proceeding as in $\S 4$ of Bourgain~\cite{bourgain}, or as at the very beginning of the proof of Proposition 3, we obtain that
$$ \delta^{2}R^{2}\Psi(x,y) \leq \sum_{1 \leq r,s \leq R} \Biggl|\sum_{\substack{n \leq x, \\ n \in \mathcal{S}(y)}} e(n(\theta_{r}-\theta_{s})) \Biggr| . $$

Now combining Theorem 1 with Dirichlet's approximation theorem, as in our previous arguments, the only way that we can possibly have $\left|\sum_{\substack{n \leq x, \\ n \in \mathcal{S}(y)}} e(n\theta) \right| \geq (\delta^{2}/2) \Psi(x,y)$ is if $\theta = a/q+\eta$ for some $q \leq (\delta^{-1} \log x)^{10}$ and some $|\eta| \leq (\delta^{-1} \log x)^{10}/x$, say. Moreover, since we have $\delta^{-1} \leq \min\{\log^{1/\epsilon}x, y^{1/C_{0}}\}$ we can apply Major Arc Estimates 1 and 2 in that case, and deduce the stronger fact that actually
\begin{eqnarray}
\sum_{\substack{n \leq x, \\ n \in \mathcal{S}(y)}} e(n(\frac{a}{q}+\eta)) & = & V(x,y;q,\eta) + O_{\epsilon}\left(\frac{\delta^{2} \Psi(x,y)}{\log x} \right) \nonumber \\
& \ll_{\epsilon} & \Psi(x,y) \frac{2^{\omega(q)} q^{1-\alpha} \log^{2}(q+1)}{\phi(q)} \frac{\log^{3}(2+|\eta x|)}{(1+|\eta x|)^{\alpha}} + \frac{\delta^{2} \Psi(x,y)}{\log x} . \nonumber
\end{eqnarray}
Here we bounded the main term $(\Psi(x,y)/x)V(x,x;q,\delta)$ in Major Arc Estimate 2 by noting it is $(\Psi(x,y)/x) \sum_{n \leq x} e(n(a/q+\eta)) + O_{\epsilon}(\delta^{2}\Psi(x,y)/\log x)$, (by Major Arc Estimate 1, again), and then summing the geometric progression.

Using the stronger bound, provided $\delta$ is small enough in terms of $\epsilon$ the only way that we can possibly have $\left|\sum_{\substack{n \leq x, \\ n \in \mathcal{S}(y)}} e(n(\theta_{r}-\theta_{s})) \right| \geq (\delta^{2}/2) \Psi(x,y)$ is actually if $\theta_{r}-\theta_{s} = a/q+\eta$ for some $q \leq \delta^{-4}$ and some $|\eta| \leq \delta^{-4}/x$, say. Thus we obtain that
\begin{eqnarray}
\delta^{2}R^{2}\Psi(x,y) & \ll_{\epsilon} & \sum_{q \leq \delta^{-4}} \sum_{(a,q)=1} \sum_{\substack{1 \leq r,s \leq R, \\ ||\theta_{r}-\theta_{s}-\frac{a}{q}|| \leq \frac{\delta^{-4}}{x}}} \Psi(x,y) \frac{2^{\omega(q)} q^{1-\alpha} \log^{2}(1+\delta^{-4}) \log^{3}(2+\delta^{-4})}{\phi(q) (1+x||\theta_{r}-\theta_{s}-a/q||)^{\alpha}} \nonumber \\
& \ll_{\epsilon} & \delta^{-8(1-\alpha)-\epsilon} \sum_{q \leq \delta^{-4}} \frac{1}{q} \sum_{(a,q)=1} \sum_{\substack{1 \leq r,s \leq R, \\ ||\theta_{r}-\theta_{s}-a/q|| \leq \delta^{-4}/x}} \frac{\Psi(x,y)}{1+x||\theta_{r}-\theta_{s}-a/q||} . \nonumber
\end{eqnarray}

The bound (\ref{largedelta}) follows swiftly from this, by applying Harmonic Analysis Result 2 with the choices $Q=\delta^{-4}$, $\Delta = \delta^{-4}/x$ and $A=1$.
\begin{flushright}
Q.E.D.
\end{flushright}

\section{Proof of Corollary 1}
The proof of Corollary 1 will be an application of the circle method, using Theorem 2, Theorem 1 and the major arc estimates in $\S 2.2$. Set $R := \log^{20}x$, and
$$ \mathfrak{M} := \bigcup_{q \leq R} \bigcup_{(a,q)=1} [a/q - R/x,a/q+R/x] . $$
In order that we may apply Theorem 1 straightforwardly, we assume throughout this section that $\log^{K}x \leq y \leq x^{1/100}$. The case of Corollary 1 where $x^{1/100} < y \leq x$ is handled by the work of La Bret\`{e}che and Granville~\cite{dlbgranville}, for example.

By the orthogonality of additive characters,
$$ \#\{(a,b,c) \in \mathcal{S}(y)^{3} : a,b,c \leq x, \; a+b=c\} = \int_{0}^{1} \Biggl(\sum_{\substack{n \leq x, \\ n \in \mathcal{S}(y)}} e(n\theta) \Biggr)^{2} \overline{\Biggl(\sum_{\substack{n \leq x, \\ n \in \mathcal{S}(y)}} e(n\theta) \Biggr)} d\theta . $$
The main term $\Psi(x,y)^{3}/2x$ in Corollary 1 comes from the integral over $\mathfrak{M}$, and the error term $O(\Psi(x,y)^{3}\log(u+1)/(x\log y))$ comes partly from $\mathfrak{M}$ and also from the rest of the integral. To extract the main term neatly, let us note that
$$ \frac{x^{3}}{2x}\left(1+ O(1/x) \right) = \#\{(a,b,c) : a,b,c \leq x, \; a+b=c\} = \int_{0}^{1} \left(\sum_{n \leq x} e(n\theta) \right)^{2} \overline{\left(\sum_{n \leq x} e(n\theta) \right)} d\theta , $$
by a trivial counting argument. We will show that
$$ \int_{[0,1]\backslash \mathfrak{M}} \Biggl|\sum_{\substack{n \leq x, \\ n \in \mathcal{S}(y)}} e(n\theta) \Biggr|^{3} d\theta \ll \frac{\Psi(x,y)^{3}}{x} \frac{1}{\log^{4}x} , \;\;\; \text{and} \;\;\; \int_{[0,1]\backslash \mathfrak{M}} \left|\sum_{n \leq x} e(n\theta) \right|^{3} d\theta \ll \frac{x^{3}}{x} \frac{1}{\log^{4}x} , $$
and we will also show that
\begin{eqnarray}\label{scalingrel}
\int_{\mathfrak{M}} \Biggl(\sum_{\substack{n \leq x, \\ n \in \mathcal{S}(y)}} e(n\theta) \Biggr)^{2} \overline{\Biggl(\sum_{\substack{n \leq x, \\ n \in \mathcal{S}(y)}} e(n\theta) \Biggr)} d\theta & = & \frac{\Psi(x,y)^{3}}{x^{3}} \int_{\mathfrak{M}} \left(\sum_{n \leq x} e(n\theta) \right)^{2} \overline{\left(\sum_{n \leq x} e(n\theta) \right)} d\theta \nonumber \\
&& + O\left(\frac{\Psi(x,y)^{3}}{x} \frac{\log(u+1)}{\log y}\right) .
\end{eqnarray}
Combining all these facts immediately yields the corollary.

\vspace{12pt}
Before the main proof, we shall demonstrate the following proposition.
\begin{prop5}
For any large $\log^{K}x \leq y \leq x^{1/100}$, and any $\theta \in [0,1]\backslash \mathfrak{M}$, we have
$$ \Biggl|\sum_{\substack{n \leq x, \\ n \in \mathcal{S}(y)}} e(n\theta) \Biggr| \ll \Psi(x,y) \frac{1}{\log^{5}x} . $$
\end{prop5}

To prove the proposition we simply distinguish two different cases. By Dirichlet's approximation theorem, for each $\theta \in [0,1]$ we have $\theta = a/q + \delta$ for some $q \leq x^{0.55}$, some $(a,q)=1$, and some $|\delta| \leq 1/qx^{0.55}$, and if $\theta \in [0,1]\backslash \mathfrak{M}$ then we either have $q > R$ or $|\delta| > R/x$. If $\log^{K}x \leq y \leq x^{1/100}$ and $x^{0.48} \leq q \leq x^{0.55}$ then, by Fouvry and Tenenbaum's~\cite{fouvryten0} minor arc estimate stated in the introduction,
$$ \Biggl|\sum_{\substack{n \leq x, \\ n \in \mathcal{S}(y)}} e(n\theta) \Biggr| \ll x(1+|\delta x|) \log^{3}x \left(\frac{\sqrt{y}}{x^{1/4}} + \frac{1}{\sqrt{q}} + \sqrt{\frac{qy}{x}} \right) \ll x \log^{3}x \left(\frac{1}{x^{0.245}} + \frac{1}{x^{0.24}} + \frac{1}{x^{0.22}} \right) , $$
on noting that $|\delta x| \leq x^{0.45}/q \leq 1$. Since $\Psi(x,y) \geq \Psi(x,\log^{K}x) = x^{1-1/K+o(1)}$, this estimate is much stronger than we need provided $K$ is large enough. The other case is where $\log^{K}x \leq y \leq x^{1/100}$, $q < x^{0.48}$, and either $q > R := \log^{20}x$ or $|\delta| > (\log^{20}x)/x$. Then it is easy to check that $q^{2}y^{3}(1+|\delta x|)^{2} \leq x/4$, and we have $1-\alpha(x,y) \leq 1/K + o(1) \leq 1/100$ (say) by the approximation (\ref{alphaapprox}), so Theorem 1 is applicable and yields that
$$ \left|\sum_{\substack{n \leq x, \\ n \in \mathcal{S}(y)}} e(n\theta) \right| \ll \frac{\Psi(x,y)}{(\log^{20}x)^{1/2 - (3/2)(1-\alpha)}} \log^{7/2}x \ll \frac{\Psi(x,y)}{(\log^{20}x)^{0.45}} \log^{7/2}x \ll \frac{\Psi(x,y)}{\log^{5}x} . $$
\begin{flushright}
Q.E.D.
\end{flushright}

Now we can deduce our claimed bounds for the integrals over $[0,1]\backslash \mathfrak{M}$. Provided $y \geq \log^{K}x$ with $K$ large enough, Theorem 2 implies that
$$ \int_{[0,1]\backslash \mathfrak{M}} \Biggl|\sum_{\substack{n \leq x, \\ n \in \mathcal{S}(y)}} e(n\theta) \Biggr|^{3} \leq \sup_{\theta \notin \mathfrak{M}} \Biggl|\sum_{\substack{n \leq x, \\ n \in \mathcal{S}(y)}} e(n\theta) \Biggr|^{0.9} \int_{0}^{1} \Biggl|\sum_{\substack{n \leq x, \\ n \in \mathcal{S}(y)}} e(n\theta) \Biggr|^{2.1} \ll \frac{\Psi(x,y)^{2.1}}{x} \sup_{\theta \notin \mathfrak{M}} \Biggl|\sum_{\substack{n \leq x, \\ n \in \mathcal{S}(y)}} e(n\theta) \Biggr|^{0.9} . $$
Then using Proposition 5, the above is all $\ll (\Psi(x,y)^{3}/x)(1/\log^{4}x)$, as claimed. The bound $\int_{[0,1]\backslash \mathfrak{M}} \left|\sum_{n \leq x} e(n\theta) \right|^{3} d\theta \ll \frac{x^{3}}{x} \frac{1}{\log^{4}x}$ for the complete sum is a trivial consequence of the pointwise bound $\sum_{n \leq x} e(n\theta) \ll 1/||\theta||$.

It only remains to show the approximate scaling relation (\ref{scalingrel}) between the integrals over $\mathfrak{M}$. But when $\theta = a/q+\delta \in \mathfrak{M}$ (so that $q, |\delta x| \leq R:= \log^{20}x$), and provided $y \geq \log^{K}x$ with $K$ large enough, Major Arc Estimates 1 and 2 imply that
\begin{eqnarray}
\sum_{\substack{n \leq x, \\ n \in \mathcal{S}(y)}} e(n\theta) &= & \frac{\Psi(x,y)}{x} V(x,x;q,\delta) + O\left(\Psi(x,y) \frac{\log(u+1)}{\log y} \frac{2^{\omega(q)} q^{1-\alpha} \log^{2}(q+1)}{\phi(q)} \frac{\log^{3}(2+|\delta x|)}{(1+|\delta x|)^{\alpha}} \right) \nonumber \\
& = & \frac{\Psi(x,y)}{x} \sum_{n \leq x} e(n\theta) + O\left(\Psi(x,y) \frac{\log(u+1)}{\log y} \frac{2^{\omega(q)} q^{1-\alpha} \log^{2}(q+1)}{\phi(q)} \frac{\log^{3}(2+|\delta x|)}{(1+|\delta x|)^{\alpha}} \right) . \nonumber
\end{eqnarray}
Here the second equality follows by using Major Arc Estimate 1 again, to compare $V(x,x;q,\delta)$ with $\sum_{n \leq x} e(n\theta)$. Moreover, we certainly always have $\sum_{n \leq x} e(n\theta) \ll \min\{x,1/||\theta||\} \ll x \frac{2^{\omega(q)} q^{1-\alpha} \log^{2}(q+1)}{\phi(q)} \frac{\log^{3}(2+|\delta x|)}{(1+|\delta x|)^{\alpha}}$ when $\theta \in \mathfrak{M}$, and therefore
\begin{eqnarray}
&& \Biggl| \int_{\mathfrak{M}} \Biggl(\sum_{\substack{n \leq x, \\ n \in \mathcal{S}(y)}} e(n\theta) \Biggr)^{2} \overline{\Biggl(\sum_{\substack{n \leq x, \\ n \in \mathcal{S}(y)}} e(n\theta) \Biggr)} d\theta - \frac{\Psi(x,y)^{3}}{x^{3}} \int_{\mathfrak{M}} \left(\sum_{n \leq x} e(n\theta) \right)^{2} \overline{\left(\sum_{n \leq x} e(n\theta) \right)} d\theta \Biggr| \nonumber \\
& \ll & \Psi(x,y)^{3} \frac{\log(u+1)}{\log y} \sum_{q \leq R} \sum_{(a,q)=1} \left( \frac{2^{\omega(q)} q^{1-\alpha} \log^{2}(q+1)}{\phi(q)} \right)^{3} \int_{-R/x}^{R/x} \left(\frac{\log^{3}(2+|\delta x|)}{(1+|\delta x|)^{\alpha}} \right)^{3} d\delta . \nonumber
\end{eqnarray}
Since $\alpha \geq 0.99$ (say), the double sum over $q$ and $a$ is $\ll \sum_{q \leq R} 1/q^{3/2} \ll 1$ and the integral over $\delta$ is $\ll \int 1/(1+|\delta x|)^{2} \ll 1/x$. This gives the desired relation (\ref{scalingrel}).
\begin{flushright}
Q.E.D.
\end{flushright}

As we noted in the Introduction, the error term $O(\log(u+1)/\log y)$ in Corollary 1 only tends to zero if $(\log y)/\log\log x \rightarrow \infty$. To obtain a negligible error term for smaller $y$ (but still satisfying $y \geq \log^{K}x$), one can replace Major Arc Estimate 2 in our analysis by Proposition 2.3 of Drappeau~\cite{drappeausommes}. There the main term is no longer $(\Psi(x,y)/x) V(x,x;q,\delta)$, but is a more complicated object whose contribution could be analysed as in the conditional work of Lagarias and Soundararajan~\cite{lagariassound0,lagariassound}. We leave this to the interested reader. The point here is that the precise main term for Corollary 1 itself involves the saddle-point $\alpha(x,y)$, and when $1-\alpha(x,y) = \frac{\log(u\log u) + O(1)}{\log y}$ becomes larger the simple version $\Psi(x,y)^{3}/2x$ that we stated loses accuracy.

\section{Proof of Corollary 2}
To prove Corollary 2, we shall combine Theorem 2, the major arc estimates in $\S 2.2$, and the following result, which is essentially proved in $\S 6$ of Green's paper~\cite{greenroth} but whose precise formulation we take from a slightly later paper of Green and Tao~\cite{greentaorest}.
\begin{trans1}[See Proposition 5.1 of Green and Tao~\cite{greentaorest}]
Let $N$ be a large prime, and let $0 < \delta \leq 1$. Let $f : \Z/N\Z \rightarrow [0,\infty)$ and $\nu : \Z/N\Z \rightarrow [0,\infty)$ be any functions such that
$$ f(n) \leq \nu(n) \; \forall n \in \Z/N\Z, \;\;\;\;\; \text{and} \;\;\;\;\; \frac{1}{N} \sum_{n \in \Z/N\Z} f(n) \geq \delta . $$
In addition, let $\eta \geq 0$ and $M > 0$ and $2 < p < 3$ be any parameters, and suppose that
$$ \left|\frac{1}{N} \sum_{n \in \Z/N\Z} \nu(n) e(\frac{an}{N}) - \textbf{1}_{a=0} \right| \leq \eta \; \forall a \in \Z/N\Z, \;\;\; \text{and} \;\;\; \sum_{a \in \Z/N\Z} \left|\frac{1}{N} \sum_{n \in \Z/N\Z} f(n) e(\frac{an}{N})\right|^{p} \leq M , $$
where $\textbf{1}$ denotes the indicator function.

Then
$$ \frac{1}{N^{2}} \sum_{n,d \in \Z/N\Z} f(n)f(n+d)f(n+2d) \geq c(\delta) - O_{\delta,M,p}(\eta) , $$
where $c(\delta) > 0$ depends on $\delta$ only.
\end{trans1}

Transference Principle 1 asserts that if the non-negative function $f$ has average value at least $\delta$, and the average of the $p$-th power of $\sum_{n \in \Z/N\Z} f(n) e(an/N)$ is suitably bounded, and if $f$ is {\em upper bounded} by a function $\nu$ whose exponential sums (except the trivial one at $a=0$) are small enough, then the function $f$ will count quite a lot of three term arithmetic progressions $n,n+d,n+2d$. Note that $f$ and $\nu$ are defined on the additive group $\Z/N\Z$, so Transference Principle 1 actually counts arithmetic progressions modulo $N$. This means that we cannot immediately apply the result to prove Corollary 2, (since the progressions it supplies may ``wrap around'' and not be genuine progressions of integers), but there is a standard and easy trick to avoid this problem, as we shall explain below.

For the benefit of an unfamiliar reader, we remark that if $\nu \equiv 1$ then one can take $\eta = 0$, and one can take $p=2.1$ (say) and $M=1$ in view of Parseval's identity. Thus Transference Principle 1 readily implies the original form of Roth's theorem, that positive density subsets of the integers contain three term arithmetic progressions. For sparse sets such as the smooth numbers (or the primes), one needs to choose $f$ and $\nu$ as functions that grow unboundedly with $N$, and so verifying the conditions of Transference Principle 1 becomes more challenging.

Recall that in Corollary 2 we are given $B \subseteq \mathcal{S}(y) \cap [1,x]$ such that $\#B \geq \beta \Psi(x,y)$. Choose any large prime $2x < N < 4x$, say, and note that if $n,n+d,n+2d \in B$ modulo $N$ then we must in fact have $n,n+d,n+2d \in B$ as integers, since there is too much space modulo $N$ for any wraparound to occur. In particular, if we slightly abuse notation, and let $B : \Z/N\Z \rightarrow \{0,1\}$ also denote the characteristic function of the set $B$ (modulo $N$), then to prove Corollary 2 it will suffice to show that
$$ \sum_{n \in \Z/N\Z} \sum_{\substack{d \in \Z/N\Z, \\ d \neq 0}} B(n)B(n+d)B(n+2d) \geq 1 . $$

We can simply set
$$ \nu(n) := \frac{N}{\Psi(N,y)} \textbf{1}_{n \leq N, n \in \mathcal{S}(y)}, \;\;\;\;\; \text{and} \;\;\;\;\; f(n) := \frac{N}{\Psi(N,y)} B(n) , $$
where $\textbf{1}$ denotes the indicator function. Thus by construction we have $\frac{1}{N} \sum_{n \in \Z/N\Z} \nu(n) = 1$. Moreover, if $a \neq 0$ and if $a/N$ {\em is not} within distance $(\log^{20}N)/N$ of a rational with denominator $\leq \log^{20}N$, then Proposition 5 from the previous section yields that
$$ \frac{1}{N} \sum_{n \in \Z/N\Z} \nu(n) e(\frac{an}{N}) = \frac{1}{\Psi(N,y)} \sum_{\substack{n \leq N, \\ n \in \mathcal{S}(y)}} e(\frac{an}{N}) \ll \frac{1}{\log^{5}N} \leq \frac{1}{\log^{5}x} . $$
If instead $a \neq 0$ and $a/N$ {\em is} within distance $(\log^{20}N)/N$ of a rational with denominator $\leq \log^{20}N$, then Major Arc Estimates 1 and 2 imply that
$$ \frac{1}{N} \sum_{n \in \Z/N\Z} \nu(n) e(\frac{an}{N}) = \frac{1}{\Psi(N,y)} \sum_{\substack{n \leq N, \\ n \in \mathcal{S}(y)}} e(\frac{an}{N}) = \frac{1}{N} \sum_{n \leq N} e(\frac{an}{N}) + O\left(\frac{\log(u+1)}{\log y}\right) \ll \frac{\log\log x}{\log y} , $$
since the final complete sum over $n$ is identically zero.
 
We obviously always have $f(n) \leq \nu(n)$, and we also note that
$$ \frac{1}{N} \sum_{n \in \Z/N\Z} f(n) = \frac{\#B}{\Psi(N,y)} \geq \frac{\beta \Psi(x,y)}{\Psi(N,y)} \gg \beta , $$
using Smooth Numbers Result 2 and the fact that $N < 4x$. Moreover, we have
\begin{eqnarray}
\sum_{a \in \Z/N\Z} \left|\frac{1}{N} \sum_{n \in \Z/N\Z} f(n) e(an/N)\right|^{2.1} & = & \frac{1}{\Psi(N,y)^{2.1}} \sum_{a \in \Z/N\Z} \left|\sum_{n \in B} e(an/N) \right|^{2.1} \nonumber \\
& \ll & \frac{N}{\Psi(N,y)^{2.1}} \int_{0}^{1} \left|\sum_{n \in B} e(n\theta) \right|^{2.1}d\theta \nonumber \\
& \ll & 1 . \nonumber
\end{eqnarray}
Here the first inequality is a general result of Marcinkiewicz and Zygmund, which the reader may find as Lemma 6.5 of Green~\cite{greenroth} (and which is hopefully intuitively plausible, since we might imagine that often $\left|\sum_{n \in B} e(an/N) \right|^{2.1} \approx N \int_{a/N-1/2N}^{a/N+1/2N} \left|\sum_{n \in B} e(n\theta) \right|^{2.1}d\theta$), and the final inequality follows from Theorem 2.

In summary, we have shown that all the hypotheses of Transference Principle 1 are satisfied, with $\delta = \beta/1000$ (say) and with $p=2.1$, $M=O(1)$, and $\eta = O((\log\log x)/\log y)$, so we conclude that
\begin{eqnarray}
\sum_{n \in \Z/N\Z} \sum_{d \in \Z/N\Z} B(n)B(n+d)B(n+2d) & = & \left(\frac{\Psi(N,y)}{N}\right)^{3} \sum_{n,d \in \Z/N\Z} f(n)f(n+d)f(n+2d) \nonumber \\
& \geq & \frac{\Psi(N,y)^{3}}{N} (c(\delta)-O_{\delta}(\frac{\log\log x}{\log y})) . \nonumber
\end{eqnarray}
Provided $x$ and $(\log y)/\log\log x$ are large enough in terms of $\beta$, the right hand side will be $\gg_{\beta} \Psi(x,y)^{3}/x$, whilst the trivial contribution to the left hand side (from $d=0$) is equal to $\#B \leq \Psi(x,y)$. In particular, there will be at least one (and in fact very many) non-trivial three term progressions in $B$.
\begin{flushright}
Q.E.D.
\end{flushright}

\appendix

\section{Proofs of the Major Arc Estimates}
In this appendix we shall prove the Major Arc Estimates that we stated in $\S 2.2$. In all cases, a large portion of the result can be imported immediately from the existing literature (primarily the paper of Drappeau~\cite{drappeausommes}), leaving only a few bad terms or awkward cases to be handled.

\subsection{Proof of Major Arc Estimate 1}
Under the conditions of Major Arc Estimate 1, if we apply Proposition 2.2 of Drappeau~\cite{drappeausommes} with the choice $Q := \min\{y^{c_{2}},e^{c_{2}\sqrt{\log x}}\}$ we obtain (after suitably relabelling the constant $c_{2}$) that
$$ \sum_{\substack{n \leq x, \\ n \in \mathcal{S}(y)}} e(n(a/q+\delta)) = V(x,y;q,\delta) + O\left(|W(x,y;q,\delta)| + \Psi(x,y)(1+|\delta x|)( y^{-c_{2}} + e^{-2c_{2}\sqrt{\log x}})\right) . $$
We will explain the new term $W(x,y;q,\delta)$ in a moment, but note immediately that for any $A > 0$ we have
$$ e^{-2c_{2}\sqrt{\log x}} = O_{A}(\frac{e^{-c_{2}\sqrt{\log x}}}{\log^{A}x}) = O_{A}(y^{-c_{2}} + \frac{e^{-c_{2}u/\log^{2}(u+1)}}{\log^{A}x}) , $$
since $u:=(\log x)/\log y$ and therefore we always have $e^{-\sqrt{\log x}} \leq \max\{y^{-1},e^{-u}\}$. This contribution is acceptable for Major Arc Estimate 1, so to prove the estimate it only remains to show that, for any $A > 0$,
\begin{equation}\label{needforw}
|W(x,y;q,\delta)| \ll_{A} \Psi(x,y)(1+|\delta x|)(y^{-c_{2}} + \frac{e^{-c_{2}u/\log^{2}(u+1)}}{\log^{A}x}) .
\end{equation}

As explained by Drappeau~\cite{drappeausommes}, the term $W(x,y;q,\delta)$ may be omitted unless there exists a real primitive Dirichlet character $\chi_{\text{bad}}$, with conductor $r_{\text{bad}} \leq Q = \min\{y^{c_{2}},e^{c_{2}\sqrt{\log x}}\}$, such that $r_{\text{bad}}$ divides $q$ and such that the Dirichlet $L$-function $L(s,\chi_{\text{bad}})$ has a real zero that is $\geq 1-c/\log Q$ (for a certain absolute constant $c>0$ which ensures there is at most one such primitive character). If such an exceptional character actually exists, and if we let $\tau(\cdot)$ denote the Gauss sum, and let $\chi_{\text{bad}}^{(d)}$ denote the character to modulus $dr_{\text{bad}}$ that is induced by $\chi_{\text{bad}} = \chi_{\text{bad}}^{(1)}$, then $W(x,y;q,\delta)$ is defined by
$$ W(x,y;q,\delta) := \frac{\tau(\chi_{\text{bad}})}{\phi(r_{\text{bad}})} \sum_{d|(q/r_{\text{bad}})} \frac{\mu(d)\chi_{\text{bad}}(d)}{\phi(d)} \sum_{\substack{n \leq x r_{\text{bad}}d/q, \\ n \in \mathcal{S}(y)}} e(\delta nq/dr_{\text{bad}}) \chi_{\text{bad}}^{(d)}(n) . $$
Using the standard bound $|\tau(\chi_{\text{bad}})| \leq \sqrt{r_{\text{bad}}}$ (as in e.g. chapter 9 of Davenport~\cite{davenport}), we quickly obtain that
$$ |W(x,y;q,\delta)| \ll \frac{\sqrt{r_{\text{bad}}} \log x}{\phi(r_{\text{bad}})} \max_{d|(q/r_{\text{bad}})} \left| \sum_{\substack{n \leq x r_{\text{bad}}d/q, \\ n \in \mathcal{S}(y)}} e(\delta nq/dr_{\text{bad}}) \chi_{\text{bad}}^{(d)}(n) \right| . $$
And now we claim that to prove Major Arc Estimate 1, it will suffice to show that
$$ \max_{d|(q/r_{\text{bad}})} \left| \sum_{\substack{n \leq x r_{\text{bad}}d/q, \\ n \in \mathcal{S}(y)}} e(\delta nq/dr_{\text{bad}}) \chi_{\text{bad}}^{(d)}(n) \right| \ll \Psi(x,y) (1+|\delta x|) \log^{2}x ( y^{-2c_{2}} + e^{-2c_{2}u/\log^{2}(u+1)}) . $$
For if $\log^{c_{1}}x \leq y \leq e^{\sqrt{\log x}}$, say, then we have $y^{-c_{2}} \ll 1/\log^{3}x$ and $e^{-c_{2}u/\log^{2}(u+1)} \ll_{A} 1/\log^{A+3}x$, and so our desired bound (\ref{needforw}) will follow. And if instead $e^{\sqrt{\log x}} < y \leq x$, so that $\log Q= c_{2}\sqrt{\log x}$, then any exceptional conductor will satisfy $r_{\text{bad}} \gg_{A} \log^{A}x$ for any fixed $A > 0$, by Siegel's theorem (as in e.g. chapter 21 of Davenport~\cite{davenport}). Therefore we will get the bound (\ref{needforw}) because of the prefactor $\frac{\sqrt{r_{\text{bad}}}}{\phi(r_{\text{bad}})} \ll_{A} 1/\log^{A}x$.

Next, note that if $q/r_{\text{bad}}d \geq \min\{y^{3c_{2}},e^{3c_{2}u/\log^{2}(u+1)}\}$ then the trivial bound
$$ \left| \sum_{\substack{n \leq x r_{\text{bad}}d/q, \\ n \in \mathcal{S}(y)}} e(\delta nq/dr_{\text{bad}}) \chi_{\text{bad}}^{(d)}(n) \right| \leq \Psi(x r_{\text{bad}}d/q,y) \ll \frac{\Psi(x,y)}{\min\{y^{3c_{2}},e^{3c_{2}u/\log^{2}(u+1)}\}^{\alpha(x,y)}} $$
is sufficient, (here the second inequality is Smooth Numbers Result 2), since the saddle-point $\alpha(x,y) \geq \alpha(x,\log^{c_{1}}x) \geq 2/3$ by the approximation (\ref{alphaapprox}).

If instead $q/r_{\text{bad}}d < \min\{y^{3c_{2}},e^{3c_{2}u/\log^{2}(u+1)}\}$ then we may split into subsums each of length at most 
$M := \max\{x/y^{5c_{2}},x/e^{5c_{2}u/\log^{2}(u+1)}\}$, say, and for any subsum we will have
\begin{eqnarray}
\sum_{\substack{x' < n \leq x'+M, \\ n \in \mathcal{S}(y)}} e(\frac{\delta nq}{dr_{\text{bad}}}) \chi_{\text{bad}}^{(d)}(n) & = & e(\frac{\delta x'q}{dr_{\text{bad}}}) \sum_{\substack{x' < n \leq x'+M, \\ n \in \mathcal{S}(y)}} \chi_{\text{bad}}^{(d)}(n) + O\Biggl(\frac{Mq|\delta|}{dr_{\text{bad}}} \sum_{\substack{x' < n \leq x'+M, \\ n \in \mathcal{S}(y)}} 1 \Biggr) \nonumber \\
& = &  e(\frac{\delta x'q}{dr_{\text{bad}}}) \sum_{\substack{x' < n \leq x'+M, \\ n \in \mathcal{S}(y)}} \chi_{\text{bad}}^{(d)}(n) + \nonumber \\
&& + O\Biggl(|\delta x|(y^{-2c_{2}} + e^{-2c_{2}u/\log^{2}(u+1)}) \sum_{\substack{x' < n \leq x'+M, \\ n \in \mathcal{S}(y)}} 1 \Biggr) . \nonumber
\end{eqnarray}
Here the ``big Oh'' term is acceptably small.

It only remains to bound the untwisted character sums $\sum_{\substack{x' < n \leq x'+M, \\ n \in \mathcal{S}(y)}} \chi_{\text{bad}}^{(d)}(n)$. However, this has actually already been done in previous work of the author~\cite{harpersmoothbv}. Indeed, if $\log^{c_{1}}x \leq y \leq x^{1/(\log\log x)^{2}}$ (say) then the exceptional characters argument at the end of $\S 3.3$ of that paper immediately yields that
\begin{eqnarray}
\sum_{\substack{x' < n \leq x'+M, \\ n \in \mathcal{S}(y)}} \chi_{\text{bad}}^{(d)}(n) = \sum_{\substack{n \leq x'+M, \\ n \in \mathcal{S}(y)}} \chi_{\text{bad}}^{(d)}(n) - \sum_{\substack{n \leq x', \\ n \in \mathcal{S}(y)}} \chi_{\text{bad}}^{(d)}(n) & \ll & \Psi(x'+M,y)\log^{2}x \left(y^{-c} + e^{-c u/\log^{2}(u+1)}\right) \nonumber \\
& \ll & \Psi(x,y)\log^{2}x \left(y^{-c} + e^{-c u/\log^{2}(u+1)}\right) , \nonumber
\end{eqnarray}
for a small absolute constant $c > 0$. This is clearly sufficient, on summing over the $O(\min\{y^{5c_{2}},e^{5c_{2}u/\log^{2}(u+1)}\})$ values of $x'$, provided our constant $c_{2}$ is small enough in terms of $c$. If instead $x^{1/(\log\log x)^{2}} < y \leq x$, then the argument from $\S 3.3$ of \cite{harpersmoothbv} implies that
$$ \sum_{\substack{x' < n \leq x'+M, \\ n \in \mathcal{S}(y)}} \chi_{\text{bad}}^{(d)}(n) \ll \left(\sum_{b|d} \frac{1}{b^{\alpha}} \right) \Psi(x'+M,y) \frac{\log(r_{\text{bad}})}{\log x} e^{-c u/\log^{2}(u+1)} \ll (\sum_{b|d} \frac{1}{b^{\alpha}} ) \Psi(x,y) e^{-c u/\log^{2}(u+1)} . $$
However, since $d \leq q \leq y^{c_{2}}$ for some small $c_{2}$, and since $\alpha = \alpha(x,y) = 1 + O(\log(u+1)/\log y)$ by the approximation (\ref{alphaapprox}), the sum over $b$ is easily seen to be $\ll u \log x \ll \log^{2}x$, so again our bound is sufficient.
\begin{flushright}
Q.E.D.
\end{flushright}

\subsection{Proof of Major Arc Estimate 2}
Suppose first that $q=1$, in which case the estimate we are required to prove is that
\begin{equation}\label{1case}
\sum_{\substack{n \leq x, \\ n \in \mathcal{S}(y)}} e(n\delta) = \frac{\Psi(x,y)}{x} \sum_{n \leq x} e(n\delta) + O\left(\Psi(x,y) \frac{\log(u+1)}{\log y} \frac{\log^{3}(2+|\delta x|)}{(1+|\delta x|)^{\alpha}} \right) ,
\end{equation}
whenever $\log^{c_{1}}x \leq y \leq x$ and $|\delta| \leq \min\{y^{c_{2}}/x, e^{c_{2}\log^{1/4}x}/x\}$. The result is trivial if $\delta = 0$, so suppose throughout that $|\delta| > 0$. Also, let us recall throughout that the saddle-point $\alpha=\alpha(x,y)$ satisfies $\alpha = 1 - O(\frac{\log(u+1)}{\log y})$, by the approximation (\ref{alphaapprox}), where $u := (\log x)/\log y$.

If $\log^{c_{1}}x \leq y \leq x^{1/(\log\log x)^{4}}$ and $|\delta| \leq \min\{y^{c_{2}}/x, e^{c_{2}(\log^{3/5}x)/(\log\log x)^{1/5}}/x\}$, for suitable values of $c_{1},c_{2}$, then Proposition 2.3 of Drappeau~\cite{drappeausommes} implies that
$$ \sum_{\substack{n \leq x, \\ n \in \mathcal{S}(y)}} e(n\delta) = \alpha \Psi(x,y) \int_{0}^{1} e(\delta xt) t^{\alpha-1} dt + O\left(\Psi(x,y) \frac{1}{u} \frac{\log^{3}(2+|\delta x|)}{(1+|\delta x|)^{\alpha}} \right) . $$
In addition, integration by parts shows that
\begin{eqnarray}
\int_{0}^{1} e(\delta xt) t^{\alpha-1} dt & = & \left\lbrack \frac{e(\delta xt) - 1}{2\pi i \delta x} t^{\alpha-1} \right\rbrack_{0}^{1} - (\alpha - 1) \int_{0}^{1} \frac{e(\delta xt) - 1}{2\pi i \delta x} t^{\alpha-2} dt \nonumber \\
& = & \frac{e(\delta x) - 1}{2\pi i \delta x} + O\left((1-\alpha)\int_{0}^{1} \min\{t,\frac{1}{|\delta x|}\} t^{\alpha-2}dt \right) . \nonumber
\end{eqnarray}
Here the integral is clearly $O(1)$ if $|\delta x| \leq 1$, whilst if $1 < |\delta x|$ we split the integral at $t=1/|\delta x|$ and find it is
$$ \ll \frac{1}{|\delta x|^{\alpha}} + \frac{|\delta x|^{1-\alpha}-1}{(1-\alpha)|\delta x|} \ll \frac{\log(2+|\delta x|)}{(1+|\delta x|)^{\alpha}} , $$
by distinguishing cases according as $(1-\alpha)\log(2+|\delta x|)$ is large or small.

Putting everything together, and also using the bounds $1-\alpha = O(\frac{\log(u+1)}{\log y})$ and $\left|\frac{e(\delta x) - 1}{2\pi i \delta x}\right| \ll 1/(1+|\delta x|)$, it follows that
\begin{eqnarray}
\sum_{\substack{n \leq x, \\ n \in \mathcal{S}(y)}} e(n\delta) & = & \alpha \Psi(x,y) \frac{e(\delta x) - 1}{2\pi i \delta x} + O\left(\Psi(x,y) \left(\frac{1}{u} + \frac{\log(u+1)}{\log y}\right) \frac{\log^{3}(2+|\delta x|)}{(1+|\delta x|)^{\alpha}} \right) \nonumber \\
& = & \Psi(x,y) \frac{e(\delta x) - 1}{2\pi i \delta x} + O\left(\Psi(x,y) \left(\frac{1}{u} + \frac{\log(u+1)}{\log y}\right) \frac{\log^{3}(2+|\delta x|)}{(1+|\delta x|)^{\alpha}} \right) . \nonumber
\end{eqnarray}
By summing the geometric progression we find $\frac{e(\delta x) - 1}{2\pi i \delta x} = (1/x) \sum_{n \leq x} e(n\delta) + O(1/x)$. Thus we obtain the claimed estimate (\ref{1case}) provided that $\log^{c_{1}}x \leq y \leq e^{\sqrt{\log x}}$ (say), so that $1/u \ll \log(u+1)/\log y$ in the error term.

We still need to handle the range where $e^{\sqrt{\log x}} < y \leq x$ and $|\delta| \leq \min\{y^{c_{2}}/x, e^{c_{2}\log^{1/4}x}/x\} = e^{c_{2}\log^{1/4}x}/x $. To do this we apply Proposition 1 of La Bret\`{e}che~\cite{dlbAA}, which implies that
$$ \sum_{\substack{n \leq x, \\ n \in \mathcal{S}(y)}} e(n\delta) = \rho(u) \sum_{n \leq x} e(n\delta) + O\left(\Psi(x,y) \frac{\log(u+1)}{\log y} \frac{\log(2+|\delta x|)}{1+|\delta x|} \right) $$
provided $|\delta| \leq e^{c\sqrt{\log y}}/x$, for a certain small constant $c > 0$. Here $\rho(u)$ is the Dickman function, so as described in $\S 2.1$ the main term $\rho(u) \sum_{n \leq x} e(n\delta)$ is indeed
$$ \frac{\Psi(x,y)}{x} \left(1+O\left(\frac{\log(u+1)}{\log y}\right)\right) \sum_{n \leq x} e(n\delta) = \frac{\Psi(x,y)}{x} \sum_{n \leq x} e(n\delta) + O\left(\Psi(x,y) \frac{\log(u+1)}{\log y} \frac{1}{1+|\delta x|} \right) , $$
which establishes (\ref{1case}). Actually Proposition 1 of La Bret\`{e}che~\cite{dlbAA} is restricted to the range $e^{\log^{2/3+\epsilon}x} \leq y \leq x$, but almost all of the proof works under the much weaker condition that $y \geq e^{(\log\log x)^{5/3+\epsilon}}$. The assumption that $y \geq e^{\log^{2/3+\epsilon}x}$ is only made to ensure that $\Psi(x/e^{c\sqrt{\log y}},y) \ll \Psi(x,y)/e^{c\sqrt{\log y}}$, and this is certainly also true for all $e^{\sqrt{\log x}} < y \leq x$ (using e.g. Smooth Numbers Result 2 and the approximation (\ref{alphaapprox})).

It remains for us to prove Major Arc Estimate 2 in the case where $2 \leq q \leq y^{1/4}$. However, this will be a quick deduction from (\ref{1case}). For by definition we see $V(x,y;q,\delta)$ is
$$ := \sum_{\substack{n \leq x, \\ n \in \mathcal{S}(y)}} \frac{\mu(q/(q,n))}{\phi(q/(q,n))} e(n\delta) = \sum_{h|q} \frac{\mu(q/h)}{\phi(q/h)} \sum_{\substack{m \leq x/h, \\ m \in \mathcal{S}(y), \\ (m,q/h)=1}} e(mh\delta) = \sum_{h|q} \frac{\mu(q/h)}{\phi(q/h)} \sum_{d|(q/h)} \mu(d) \sum_{\substack{n \leq x/dh, \\ n \in \mathcal{S}(y)}} e(ndh\delta) , $$
where the final equality uses the fact that $\textbf{1}_{(m,q/h)=1} = \sum_{d|(m,q/h)} \mu(d)$. And we can use our result (\ref{1case}) to estimate all the inner sums, deducing that $V(x,y;q,\delta)$ is equal to
$$ \sum_{h|q} \frac{\mu(q/h)}{\phi(q/h)} \sum_{d|(q/h)} \mu(d) \left(\frac{\Psi(\frac{x}{dh},y)}{x/dh} \sum_{n \leq x/dh} e(ndh\delta) + O\left(\Psi(\frac{x}{dh},y) \frac{\log(u+1)}{\log y} \frac{\log^{3}(2+|\delta x|)}{(1+|\delta x|)^{\alpha}} \right) \right) . $$
Using Smooth Numbers Result 2 in the form $\Psi(x/dh,y) \ll \Psi(x,y)/(dh)^{\alpha} \leq \Psi(x,y)q^{1-\alpha}/dh$, we see the contribution from all the ``big Oh'' error terms is
\begin{eqnarray}
& \ll & \Psi(x,y) \frac{\log(u+1)}{\log y} \frac{\log^{3}(2+|\delta x|)}{(1+|\delta x|)^{\alpha}} q^{1-\alpha} \sum_{h|q} \frac{|\mu(q/h)|}{\phi(q/h)h} \sum_{d|(q/h)} \frac{1}{d} \nonumber \\
& \ll & \Psi(x,y) \frac{\log(u+1)}{\log y} \frac{\log^{3}(2+|\delta x|)}{(1+|\delta x|)^{\alpha}} \frac{2^{\omega(q)} q^{1-\alpha} \log(q+1)}{\phi(q)} , \nonumber
\end{eqnarray}
which is more than acceptable for Major Arc Estimate 2.

Finally, it follows from Th\'{e}or\`{e}me 2.4(ii) (last part) of La Bret\`{e}che and Tenenbaum~\cite{dlbten}, and then the approximation $\alpha = 1 - O(\log(u+1)/\log y)$, that $\Psi(x/dh,y)$ is equal to
\begin{eqnarray}
(1+O\left(\frac{\log(u+1)}{\log y} + \frac{\log(dh+1)}{\log x}\right)) \frac{\Psi(x,y)}{(dh)^{\alpha}} & = & \frac{\Psi(x,y)}{dh} (1 + O(\frac{(dh)^{1-\alpha} \log(dh+1) \log(u+1)}{\log y})) \nonumber \\
& = & \frac{\Psi(x,y)}{dh} (1 + O(\frac{q^{1-\alpha} \log(q+1) \log(u+1)}{\log y})) . \nonumber
\end{eqnarray}
Remembering that $|\sum_{n \leq x/dh} e(ndh\delta)| \ll x/(dh(1+|\delta x|))$, the contribution from all these ``big Oh'' terms to $V(x,y;q,\delta)$ may be bounded exactly as we did above. This only leaves the main term contribution, which is
$$ \sum_{h|q} \frac{\mu(q/h)}{\phi(q/h)} \sum_{d|(q/h)} \mu(d) \frac{\Psi(x,y)}{x} \sum_{n \leq x/dh} e(ndh\delta) = \frac{\Psi(x,y)}{x} V(x,x;q,\delta) , $$
exactly as claimed in Major Arc Estimate 2.
\begin{flushright}
Q.E.D.
\end{flushright}

\vspace{12pt}
\noindent {\em Acknowledgements.} The author would like to thank R\'{e}gis de la Bret\`{e}che and Sary Drappeau for their detailed comments and suggestions about a draft of this paper. He also thanks Junsoo Ha for explaining his results, and Andrew Granville, Ben Green, and K. Soundararajan for their comments and encouragement.

\end{document}